\documentclass[11pt]{amsart}
\usepackage{amssymb}
\usepackage{color}
\usepackage{verbatim} 
\newtheorem{Theorem}[equation]{Theorem}

\newtheorem{cor}[equation]{Corollary}
\newtheorem{Lemma}[equation]{Lemma}

\newtheorem{rem}[equation]{Remark}

\numberwithin{equation}{section}

\address{%
}

\def\a1s{a_1,\cdots, a_s}
\def\a{\alpha}

\def\aa{\mathcal A}

\def\fa{\mathfrak{a}}

\def\andd{\quad\hbox{and}\quad}

\def\fb{\frak{b}}
\def\b{\beta}

\def\bb{\mathcal{B}}

\def\bl4{B_{\ell\geq4}}
\def\cc{{\mathcal C}}

\def\d{\delta}

\def\dd{\mathcal D}

\def\bbbf{\mathbb{F}}

\def\gg{{\mathcal G}}

\def\fg{\mathfrak{g}}

\def\heart{\hbox{\tiny$\heartsuit$}}

\def\hh{{\mathcal H}}

\def\fh{\mathfrak{h}}

\def\ii{\mathcal{I}}
\def\jj{\mathcal{J}}

\def\kk{\mathcal{K}}

\def\lam{\lambda}
\def\Lam{\Lambda}
\def\LL{\mathcal{L}}

\def\ep{\epsilon}
\def\fm{(\cdot,\cdot)}

\def\m{\mathcal{M}}
\def\bbbn{\mathbb{N}}

\def\1k{\frac{1}{k}}
\def\op{\oplus}
\def\ot{\otimes}

\def\la{\langle}
\def\ra{\rangle}

\def\sub{\subseteq}

\def\pf{\noindent{\bf Proof. }}

\def\i{{\mathcal I}}

\def\kk{\mathcal K}

\def\ss{\mathcal{S}}

\def\T{{\mathcal T}}

\def\v{{\mathcal V}}

\def\z{{\mathcal Z}}

\def\ss{\mathcal{S}}
 
\def\1il{1\leq i\leq\ell}

\begin{document}

\centerline{\bf Central Extensions  of  Root Graded Lie Algebras}

\vspace{1cm}
\centerline{Malihe Yousofzadeh\footnote{This research was in part
supported by a grant from IPM (No. 90170031).}}



\vspace{1cm}
\parbox{5in}{}
\textbf{Abstract.} We study the central extensions of Lie algebras graded by an irreducible locally finite root system.\\
\emph{\small Mathematics Subject Classification:} 17B70, 17B60,\\
\emph{\small Key Words:} Root graded Lie algebras, Central Extensions.
\setcounter{section}{-1}
\section{introduction}
In 1992, S. Berman and R. Moody \cite{BM} introduced the notion of
a  {\it Lie algebra graded by an irreducible  reduced finite  root
system.} 
This definition was generalized by E. Neher \cite{N} in 1996
to Lie algebras graded by a  reduced  locally finite root system and also by  B. Allison, G. Benkart and Y. Gao
\cite{ABG2} in 2002 to Lie algebras   graded by an irreducible
finite root system of type $BC.$ Finally this definition was generalized to  Lie algebras graded by a locally finite root system (not necessarily reduced) in \cite{N} and \cite{Y}. 
A complete description of Lie algebras graded by an irreducible finite root system is given in \cite{BM}, \cite{BZ}, \cite{ABG2} and \cite{BS}. In \cite{N}, E. Neher  realizes Lie algebras  graded by a reduced locally finite root system other than root systems of types $F_4,$ $G_2$ and $E_8$ as  central
extensions  of  Tits-Kantor-Koecher  algebras of certain Jordan
pairs. In \cite{Y}, the author studies Lie algebras graded by an infinite irreducible locally finite root system  (not necessarily reduced) and gives a complete description of these Lie algebras.

Central extensions play a very important role in the theory of Lie algebras. Central extensions of Lie algebras graded by an irreducible finite root system is given in \cite{ABG1}, \cite{ABG2} and \cite{BS}. The universal central extension of Lie algebras graded by a reduced locally finite root system is studied by A. Walte in her Ph.D. thesis  \cite{W} in 2010. In 2011, E. Neher and  J.  Sun prove that the universal central extension of a direct limit of  a class $\{\LL_i\mid i\in I\}$ of perfect Lie superalgebras coincides with the direct limit of universal central extension of $\LL_i$'s. As a by-product, they determine the universal central extension of Lie algebras graded by an irreducible reduced locally finite root system. Here in this work we study the central extension of a Lie algebra graded by an irreducible  locally finite root system  (not necessarily reduced). According to \cite{Y}, if $X$ is the type of an irreducible locally finite root system, for  a specific quadruple $\mathfrak{q}$ called a {\it coordinate quadruple} of type $X,$ one can associate an algebra $\fb(\mathfrak{q})$ and a Lie algebra $\{\fb(\mathfrak{q}),\fb(\mathfrak{q})\}.$ The structure of  Lie algebras graded by an irreducible locally finite root system $R$ of type $X$  just depends on coordinate  quadruples $\mathfrak{q}$ of  type $X$  and certain subspaces $\kk$ of $\{\fb(\mathfrak{q}),\fb(\mathfrak{q})\}$ said to satisfies  the uniform property on $\fb(\mathfrak{q}).$    In fact corresponding to a coordinate quadruple $\mathfrak{q}$ of type $X$ and a subspace $\kk$ of $\{\fb(\mathfrak{q}),\fb(\mathfrak{q})\} $ satisfying the uniform property on $\fb(\mathfrak{q}),$ the author associates a Lie algebra $\LL(\mathfrak{q},\kk)$ and shows that it is a Lie algebra graded by the irreducible locally finite root system $R$ of type $X.$ Conversely she proves that any $R-$graded Lie algebras is isomorphic to such a Lie algebra. In this work we study the central extensions of root graded Lie algebras. We prove that a perfect central extension of a Lie algebra graded by  an irreducible locally finite root system $R$  is a Lie algebra graded by  the  same root system $R$ with the same coordinate quadruple. Moreover, we prove that the universal central extension of a Lie algebra $\LL=\LL(\mathfrak{q},\kk)$ graded by an irreducible locally finite root system $R$ is $\LL(\mathfrak{q}, \{0\}).$ 

\section{Preliminary}

By a {\it star algebra} $(\mathfrak{A},\star),$ we mean an algebra $\mathfrak{A}$ together
with   a self-inverting antiautomorphism $\star$ which is referred  to as an {\it involution}.

We call a quadruple $(\fa,*,\cc,f),$ a {\it coordinate  quadruple} if one of the followings holds:
\begin{itemize}
\item  (Type $A$) $\fa$ is a unital associative algebra, $*=id_\fa,$ $\cc=\{0\}$ and $f:\cc\times\cc\longrightarrow \fa$ is the  zero map.
\smallskip

\item    (Type $B$) $\fa=\aa\op\bb$  where $\aa$ is a unital
commutative associative algebra and $\bb$ is a unital associative
$\aa$-module equipped with a symmetric bilinear form and $\fa$ is
the corresponding Clifford  Jordan algebra, $*$ is a linear
transformation fixing the elements of $\aa$ and skew fixing the
elements of $\bb,$ $\cc=\{0\}$ and $f:\cc\times\cc\longrightarrow
\fa$ is the zero map.
\smallskip

\item  (Type $C$) $\fa$ is a unital associative algebra, $*$ is an involution on $\fa,$ $\cc=\{0\}$ and $f:\cc\times\cc\longrightarrow \fa$ is the zero map.
\smallskip

\item  (Type $D$) $\fa$ is a unital  commutative associative algebra $*=id_\fa,$ $\cc=\{0\}$ and $f:\cc\times\cc\longrightarrow \fa$ is the zero map.
\smallskip

\item  (Type $BC$) $\fa$ is a unital associative algebra, $*$ is
an involution on $\fa,$ $\cc$ is a unital associative $\fa$-module
and $f:\cc\times\cc\longrightarrow \fa$ is a skew-hermitian form.

\end{itemize}

Suppose that $\mathfrak{q}:=(\fa,*,\cc,f)$ is a coordinate quadruple. Denote by  $\aa$ and $\bb,$ the fixed and the skew fixed points of  $\fa$ under $*,$ respectively.
 Set
$\fb:=\fb(\fa,*,\cc,f):=\fa\op\cc$ and define
\begin{equation}\label{probinbc-n}
\begin{array}{c}\cdot:\fb\times\fb\longrightarrow \fb\\
(\a_1+c_1,\a_2+c_2)\mapsto(\a_1\cdot \a_2)+f(c_1,c_2)+\a_1\cdot
c_2+\a_2^*\cdot c_1,
\end{array}
\end{equation}
for $\a_1,\a_2\in\fa$ and $c_1,c_2\in\cc.$
Also for $\b,\b'\in\fb,$
set
\begin{equation}\label{end9-n}
\b\circ\b':=\b\cdot\b'+\b'\cdot\b\andd [\b,\b']:=\b\cdot \b'-\b'\cdot\b,
\end{equation}
and for $c,c'\in\cc,$ define
\begin{equation}\label{diamond-heart-n}
\begin{array}{c}
\begin{array}{ll}
\diamond:\cc\times\cc\longrightarrow\aa,&
(c,c')\mapsto\frac{f(c,c')-f(c',c)}{2};\;c,c'\in\cc,
\end{array}\vspace{3mm}\\
\begin{array}{ll}
\heart:\cc\times\cc\longrightarrow\bb,&
(c,c')\mapsto\frac{f(c,c')+f(c',c)}{2};\;c,c'\in\cc.
\end{array}
\end{array}
\end{equation}

Now suppose that $\ell$ is a positive integer and for $\a,\a'\in\fa$ and $c,c'\in\cc,$ consider the following
endomorphisms
\begin{equation}\label{derivbc}
\begin{array}{l}
d_{\a,\a'}:\fb\longrightarrow\fb,\\
\b\mapsto\left\{\begin{array}{ll}
\frac{1}{\ell+1}[[\a,\a'],\b]& \hbox{$\mathfrak{q}$  is of type $A,$}\; \b\in \fb,\vspace{2mm}\\
\a'(\a\b)-\a(\a'\b)&\hbox{$\mathfrak{q}$  is of type $B,$}\;\b\in\fb,\vspace{2mm}\\
 \frac{1}{4\ell}[[\a,\a']+[\a^{*},\a'^{*}],\b]&\hbox{$\mathfrak{q}$  is of type $C$ or $BC,$}\;\;\b\in\fa,\vspace{2mm}\\
\frac{1}{4\ell}([\a,\a']+[\a^{*},\a'^{*}])\cdot
\b&\hbox{$\mathfrak{q}$  is of type $C$ or $BC,$}\;\;\b\in\cc,\vspace{2mm}\\
0& \hbox{$\mathfrak{q}$  is of type $D,$}\; \b\in\fb,\end{array}\right.\\d_{c,c'}:\fb\longrightarrow\fb,\\
\b\mapsto\left\{\begin{array}{ll}\frac{-1}{2\ell}[c\heart c',\b]& \hbox{$\mathfrak{q}$  is of type $BC,$}\;\b\in\fa,\vspace{2mm}\\
\frac{-1}{2\ell}(c\heart c')\cdot\b-\frac{1}{2}(f(\b,c')\cdot
c+f(\b,c)\cdot c')&\hbox{$\mathfrak{q}$  is of type $BC,$}\;\b\in\cc,\vspace{2mm}\\
0&\hbox{otherwise},\vspace{2mm}\end{array}\right.\\
d_{\a,c}:=d_{c,\a}:=0,\vspace{2mm}\\
d_{\a+c,\a'+c'}:=d_{\a,\a'}+d_{c,c'}.
\end{array}
\end{equation}
One can see that for $\b,\b'\in\fb,$ $d_{\b,\b'}\in Der(\fb).$
Next take  $K$ to be  a subspace of $\fb\ot \fb$ spanned by
$$\begin{array}{l}
\a\ot c,\;\;c\ot\a,\;\;a\ot b,\\
\a\ot\a'+\a'\ot\a,\;\;c\ot c'-c'\ot c,\\
(\a\cdot \a')\ot\a''+(\a''\cdot\a)\ot\a'+(\a'\cdot\a'')\ot\a,\\
f(c,c')\ot\a+( \a^*\cdot c')\ot c-(\a\cdot c)\ot c'
\end{array}$$
for $\a,\a',\a''\in\fa,$ $a\in\aa,$ $b\in\bb,$ and $c,c'\in\cc.$
Then $(\fb\ot\fb)/K$ is a Lie algebra under the following Lie
bracket
\begin{equation}\label{last4}
[(\b_1\ot \b_2)+K,(\b'_1\ot \b'_2)+K]:=((d_{\b_1,\b_2}(\b'_1)\ot \b'_2)+K)+(\b'_1\ot d_{\b_1,\b_2}
(\b'_2))+K)\end{equation}
for $\b_1,\b_2,\b'_1,\b'_2\in\fb$ (see \cite[Proposition 5.23]{ABG2} and \cite{ABG1}). We denote this Lie algebra by  $\{\fb,\fb\}$  (or $\{\fb,\fb\}$ if there is no confusion) and for
$\b_1,\b_2\in\fb,$ we  denote $(\b_1\ot\b_2)+K$ by   $\{\b_1,\b_2\}$ (or $\{\b_1,\b_2\}$  if there is no confusion). We recall the {\it full skew-dihedral homology group}  $${\rm HF}(\fb):=\{\sum_{i=1}^n\{\b_i,\b'_i\}\in\{\fb,\fb\}\mid \sum_{i=1}^nd_{\b_i,\b'_i}=0\}$$  of $\fb$  (with respect to $\ell$) from \cite{ABG2} and \cite{ABG1} and  note that it is a subset of the center of $\{\fb,\fb\}. $
For $\b_1=a_1+b_1+c_1\in\fb$ and $\b_2=a_2+b_2+c_2\in\fb$ with $a_1,a_2\in\aa,$ $b_1,b_2\in\bb$ and $c_1,c_2\in\cc,$ set \begin{equation}\label{beta*}\b_{_{\b_1,\b_2}}^*:=[a_1,a_2]+[b_1,b_2]-c_1\heart c_2;\;\;\b_1^*:=c_1,\;\;\b_2^*:=c_2.\end{equation}
We say a subset   $\kk$ of the full skew-dihedral homology group of $\fb$ satisfies the ``{\it uniform property on $\fb$}" if  for  $\b_1,\b'_1,\ldots,\b_n,\b'_n\in\fb,$ $\sum_{i=1}^n\{\b_i,\b'_i\}\in\kk$   implies that  $\sum_{i=1}^n\b^*_{\b_i,\b'_i}=0.$

\begin{rem}
\label{rem1}{\rm 
We point it out that if $\mathfrak{q}$ is a coordinate quadruple and ${\rm HF}(\fb(\mathfrak{q}))$ has a subspace satisfying the uniform property on $\fb(\mathfrak{q}),$ then $\{0\}$  also satisfies the uniform property on $\fb.$}
\end{rem}

Suppose that $I$ is a nonempty  index set and set $J:=I\uplus\bar I.$ Suppose that  $\v$ is a
vector space with a fixed basis $\{v_j\mid j\in J\}.$ One knows
that $\mathfrak{gl}(\v):=\hbox{End}(\v)$ together with
$$[\cdot,\cdot]:\mathfrak{gl}(\v)\times\mathfrak{gl}(\v)\longrightarrow
\mathfrak{gl}(\v);\; (X,Y)\mapsto XY-YX;\;\;
X,Y\in\mathfrak{gl}(\v)$$is a Lie algebra. Now for $j,k\in J,$
define
\begin{equation}\label{elementary2}e_{j,k}:\v\longrightarrow\v;\;\;
v_i\mapsto \d_{k,i}v_j,\;\;\; (i\in J),\end{equation} then
$\mathfrak{gl}(J):=\hbox{span}_\bbbf\{e_{j,k}\mid j,k\in J\}$ is a
Lie subalgebra of $\mathfrak{gl}(\v)$.
Consider the bilinear form $\fm$ on $\v$ defined
by\begin{equation}\label{form-c}(v_j,v_{\bar k}):=-(v_{\bar
k},v_j):=2\d_{j,k},\; (v_j,v_k):=0,\; (v_{\bar j},v_{\bar
k}):=0;\;\;\;  (j,k\in I),\end{equation}and set
$$\gg:=\mathfrak{sp}(I):=\{\phi\in\mathfrak{gl}(J)\mid
(\phi(v),w)=-(v,\phi(w)),\; \hbox{for all $v,w\in\v$}\}.$$ Also for a fixed subset $I_0$ of $I,$ take $\{I_\lam\mid
\lam\in\Lam\}$ to be the class of all finite subsets of $I$
containing $I_0,$ in which $\Lam$ is an index set containing $0,$ and  for each $\lam\in \Lam,$ set
\begin{equation}\label{simple-c-alg}\gg^\lam:=\gg\cap\hbox{span}\{e_{r,s}\mid r,s\in I_\lam\cup\bar
I_\lam\}.
\end{equation}
Then $\gg$
is a  locally finite  split simple Lie
subalgebra of $\mathfrak{gl}(J)$ with splitting Cartan subalgebra
$\hh:=\hbox{span}_\bbbf\{h_i:=e_{i,i}-e_{\bar i,\bar i}\mid i\in
I\}.$ Moreover, for $i,j\in I$ with $i\neq j,$ we have
$$\begin{array}{c}\gg_{\ep_i-\ep_j}=\bbbf (e_{i,j}-e_{\bar j,\bar i}),\;\gg_{\ep_i+\ep_j}=\bbbf (e_{i,\bar j}+e_{ j,\bar i}),\;\gg_{-\ep_i-\ep_j}=\bbbf (e_{\bar i,j}+e_{\bar j,i}),\\
\gg_{2\ep_i}=\bbbf e_{i,\bar i},\;\gg_{-2\ep_i}=\bbbf e_{\bar i,i}.\end{array}$$ Also for $\lam\in\Lam,$ $\gg^\lam$ is  a finite dimensional
split simple Lie subalgebra  of type $C,$ with  splitting Cartan
subalgebra $\hh^\lam:=\hh\cap\gg^\lam,$ and
 $\gg$ is the
direct union of $\{\gg^\lam\mid \lam\in\Lam\}.$

Define
$$\pi_1:\gg\longrightarrow \hbox{End}(\v);\;\pi(\phi)(v):=\phi(v);\;\;\phi\in\gg,\; v\in
\v.$$ Then $\pi_1$ is an irreducible  representation of $\gg$ in
$\v$ equipped with a weight space decomposition with respect to $\hh$ whose set of weights is $\{\pm\ep_i\mid i\in I\}$ with
$\v_{\ep_i}=\bbbf v_i$ and $\v_{-\ep_i}=\bbbf v_{\bar i}$ for
$i\in I.$ Also for \begin{equation}\label{module-s-c}\ss:=\{\phi\in \mathfrak{gl}(J)\mid
tr(\phi)=0,(\phi(v),w)=(v,\phi(w)),\; \hbox{for all
$v,w\in\v$}\},\end{equation} we have that  $$\pi_2:\gg\longrightarrow
\hbox{End}(\ss);\; \pi_2(X)(Y):=[X,Y];\;\; X\in \gg,\; Y\in \ss$$ is an
irreducible representation of $\gg$ in $\ss$  equipped with a weight space decomposition with respect to $\hh$ whose set of weights
is $\{0,\pm(\ep_i\pm\ep_j)\mid i,j\in I,\; i\neq j\}$ with
$\ss_0=\hbox{span}_\bbbf\{e_{r,r}+e_{\bar r,\bar
r}-\frac{1}{|I_\lam|}\sum_{i\in I_\lam} (e_{i,i}+e_{\bar i,\bar i})\mid
\lam\in\Lam,r\in I_\lam\},$ $\ss_{\ep_i+\ep_j}=\bbbf(e_{i,\bar
j}-e_{j,\bar i}),$ $\ss_{-\ep_i-\ep_j}=\bbbf(e_{\bar i,j}-e_{\bar
j, i})$ and $\ss_{\ep_i-\ep_j}=\bbbf(e_{i,j}+e_{\bar j,\bar i})$
($i,j\in I, i\neq j$).
Next for  $\lam\in \Lam,$ set
\begin{equation}\label{simple-c}\begin{array}{l}
\v^\lam:=\hbox{span}_\bbbf\{v_r\mid r\in
I_\lam\cup\bar
I_\lam\},\\\ss^\lam:=\ss\cap\hbox{span}_\bbbf\{e_{r,s}\mid
r,s\in I_\lam\cup\bar I_\lam\}.\end{array}\end{equation}   Then $\v^\lam$ and $\ss^\lam$ are  irreducible finite dimensional $\gg^\lam-$modules with the set of weights $(R_\lam)_{sh}$ and $\{0\}\cup (R_\lam)_{lg}$ respectively.

\begin{Theorem}[Recognition Theorem for  Type $BC$]\label{typebc}
Suppose that  $I$ is an infinite index set and  $\ell$ is an integer greater than 3. Assume  $R$ is an irreducible  locally finite root system of type $BC_I$ and  $\v$ is a vector space with a basis $\{v_i\mid i\in I\cup\bar I\}.$ Suppose that $\fm$ is a bilinear form as in (\ref{form-c}), set $\gg:=\mathfrak{sp}(I)$ and consider $\ss$ as in (\ref{module-s-c}). Fix a subset $I_0$ of $I$ of cardinality $\ell$ and  take $R_0$ to be the  full irreducible  subsystem  of $R$ of type $BC_{I_0}.$ Suppose that  $\{R_\lam\mid\lam\in\Lam\}$ is the class  of all finite irreducible full subsystems of $R$
containing $R_0,$ where $\Lam$ is an index set containing zero.
For $\lam\in\Lam,$ take $\gg^\lam$ as in Lemma \ref{simple-c-alg} and  $\v^\lam,\ss^\lam$ as in (\ref{simple-c}). Next  define $$\begin{array}{c}\mathfrak{I}_\lam:\v\longrightarrow\v\\
v_i\mapsto\left\{\begin{array}{ll}v_i& i\in I_\lam\cup \bar I_{\lam}\\
0&\hbox{otherwise}\end{array}\right.\end{array}$$ and for
$e,f\in\gg\cup\ss,$ define $$e\circ
f:=ef+fe-\frac{tr(ef)}{l}\mathfrak{I}_0.$$

(i) Suppose that $(\fa,*,\cc,f)$ is a coordinate quadruple of type $BC$ and $\aa,$ $\bb$ are  $*$-fixed and $*$-skew fixed points of  $\fa$ respectively. Set $\fb:=\fb(\fa,*,\cc,f)$ and take $[\cdot,\cdot], \circ,\heart,\diamond$ to be  as in Subsection \ref{subsect2-1}. For $\b_1,\b_2\in\fb,$ consider $d_{\b_1,\b_2}$  as in (\ref{derivbc}) and take $\b^*_{\b_1,\b_2},\b_1^*$ and $\b_2^*$ as in  (\ref{beta*}). For a subset $\kk$ of ${\rm HF}(\fb)$ satisfying the uniform property on $\fb,$ set $$\LL(\fb,\kk):=(\gg\ot\aa)\op(\ss\ot \bb)\op(\v\ot\cc)\op(\{\fb,\fb\}/\kk).$$  Then setting $\la \b,\b'\ra:=\{\b,\b'\}+\kk,$ $\b,\b'\in \fb,$ $\LL(\fb,\kk)$ together with
{\small\begin{equation}\label{probc-gen}
\begin{array}{l}
\;[x\ot a,y\ot a']=[x,y]\ot\frac{1}{2}(a\circ a')+ (x\circ y)\ot\frac{1}{2}[a,a']+tr(xy)\la a,a'\ra,\vspace{1mm}\\
\;[x\ot a,s\ot b]= (x\circ s)\ot\frac{1}{2}[a,b]+[x,s]\ot\frac{1}{2}(a\circ b)=-[s\ot b,x\ot a],\vspace{1mm}\\
\;[s\ot b,t\ot b']=[s,t]\ot\frac{1}{2}(b\circ b')+ (s\circ t)\ot\frac{1}{2}[b,b']+tr(st)\la b,b'\ra,\vspace{1mm}\\
\;[x\ot a,u\ot c]=xu\ot a\cdot c=-[u\ot c,x\ot a],\vspace{1mm}\\
\;[s\ot b,u\ot c]=su\ot b\cdot c=-[u\ot c,s\ot b],\vspace{1mm}\\
\;[u\ot c,v\ot c']=(u\circ v)\ot (c\diamond c')+ [u, v]\ot (c\heart c')+(u,v)\la c,c'\ra,\vspace{1mm}\\
\;[\la \b_1,\b_2\ra,x\ot a]=
\frac{-1}{4\ell}((x\circ
\mathfrak{I}_0)\ot[a,\b_{_{\b_1,\b_2}}^*]+[x,\mathfrak{I}_0]\ot (a\circ \b_{_{\b_1,\b_2}}^*)),\vspace{1mm}\\
\;[\la \b_1,\b_2\ra,s\ot
b]\hspace{-1mm}=\hspace{-1mm}\frac{-1}{4\ell}([s,\mathfrak{I}_0\hspace{-.5mm}]\hspace{-1mm}\ot\hspace{-.5mm}( b\circ \b_{_{\b_1,\b_2}}^*)\hspace{-1mm}+\hspace{-1mm}(s\circ
\mathfrak{I}_0)\hspace{-1mm}\ot \hspace{-.5mm}[b, \b_{_{\b_1,\b_2}}^*\hspace{-1mm}]\hspace{-1mm}+\hspace{-1mm}2tr(s\mathfrak{I}_0)\la b,\b_{_{\b_1,\b_2}}^*\hspace{-.5mm}\ra),\vspace{1mm}\\
\;[\la \b_1,\b_2\ra,v\ot
c]=\frac{1}{2\ell}\mathfrak{I}_0v\ot (\b_{_{\b_1,\b_2}}^*\cdot c)-\frac{1}{2}v\ot
(f(c,\b^*_2)\cdot \b^*_1+f(c,\b^*_1)\cdot \b^*_2)\\
\;[\la\b_1,\b_2\ra,\la\b'_1,\b'_2\ra]=\la d^\ell_{\b_1,\b_2}(\b'_1),\b'_2\ra+\la\b'_1,d^\ell_{\b_1,\b_2}(\b'_2)\ra
\end{array}
 \end{equation}}
for $x,y\in\gg,$  $s,t\in\ss,$ $u,v\in\v,$ $a,a'\in\aa,$
$b,b'\in\bb,$ $c,c'\in\cc,$ $\b_1,\b_2,\b_1',\b'_2\in\fb,$ is an
$R$-graded  Lie algebra with grading pair $(\gg,\hh)$ where $\hh$
is the splitting Cartan subalgebra of $\gg.$

$(ii)$ If $\LL$ is  an  $R$-graded Lie algebra with grading pair
$(\fg,\fh),$ then there is a coordinate quadruple $(\fa,*,\cc,f)$
of type $BC$ and a subspace $\kk$ of $\fb:=\fb(\fa,*,\cc,f)$
satisfying the uniform property on $\fb$ such that $\LL$ is
isomorphic to $\LL(\fb,\kk).$

\end{Theorem}

A Lie algebra epimorphism  $\pi:\tilde\LL\longrightarrow \LL$ from  $(\tilde\LL,[\cdot,\cdot\tilde])$ to $(\LL,[\cdot,\cdot])$  is called a {\it central extension of $\LL$} if  $C:=ker(\pi)\sub Z(\tilde\LL).$ One knows that there is a subspace $\LL'$ of $\tilde\LL$ such that $\pi(\LL')=\LL,$ $\pi|_{_{\LL'}}:\LL'\longrightarrow \LL$ is a linear isomorphism and $\tilde\LL=\LL'\op ker(\pi).$ For $x\in\tilde\LL,$ take  $x'\in\LL'$ and $x''\in ker(\pi)$ to be the image of $x$ under the projection maps of $\LL$ on $\LL'$ and $ker(\pi)$ respectively, then for $x,y\in \tilde\LL,$ $[x,y\tilde]=[x,y]'+[x,y]''.$
One can see that $(\LL',[\cdot,\cdot]')$ is a Lie algebra and $\pi|_{\LL'}:(\LL',[\cdot,\cdot]')\longrightarrow (\LL,[\cdot,\cdot])$ is a Lie algebra  isomorphism. Also  $\tau:\LL'\times\LL'\longrightarrow C$ mapping $(x,y)$ to $[x,y]''$ is a $2-$cocycle.
We  identify $\LL'$ with $\LL$ via $\pi,$ therefore we have  $\tilde\LL=\LL\op C,$ $\pi:\tilde\LL\longrightarrow \LL$ is the projection map and for $x,y\in\LL,$ $e,f\in C,$ $[x+e,y+f\tilde]=[x,y]+\tau(x,y).$ The central extension $\pi$ is called {\it perfect} if $\tilde\LL$ is a perfect Lie algebra.

\begin{Lemma}\label{uce3}
Suppose $\LL$ is a  Lie algebra and  $\tau:\LL\times\LL\longrightarrow C$ is a $2-$cocycle. Consider the corresponding  central extension  $\tilde\LL=\LL\op C$ with Lie bracket $[\cdot,\cdot\tilde]$ as above. Let  $\gg$ be a finite dimensional  simple Lie subalgebra of $\LL$ and consider $\tilde\LL$ as a $\gg-$module via the action action $$\begin{array}{l}\cdot:\gg\times\tilde\LL\longrightarrow \tilde\LL\\
(x,y)\mapsto [x,y\tilde];\;\;\; x\in\gg,\;\;y\in\tilde\LL.\end{array}$$ If $D$ is a trivial $\gg-$submodule of $\LL$ via the adjoint representation, then $D$ is a trivial $\gg-$submodule of $\tilde\LL,$ in particular $\tau(\gg,D)=\{0\}.$
\end{Lemma}
\pf Consider the $\gg-$submodule $D\op\tau(\gg,D)$ of $\tilde\LL.$ If $d_1,\ldots, d_n\in D$ and $r_1,\ldots,r_n\in\tau(\gg,D),$ then $\{d_1+r_1,\ldots,d_n+r_n\}$ is a subset of $\hbox{span}_\bbbf\{d_1,\ldots,d_n\}+\sum_{i=1}^n\tau(\gg,d_i)+\hbox{span}_\bbbf\{r_1,\ldots,r_n\}$ which is a finite dimensional $\gg-$submodule of $D\op\tau(\gg,D).$ This means that $D\op\tau(\gg,D)$ is a locally finite $\gg-$module and so it is completely reducible as $\gg$ is a finite dimensional simple Lie algebra. Next we note that $\tau(\gg,D)$ is a trivial $\gg-$submodule of $\dd\op\tau(\gg,\dd),$ so   there is a submodule $\dot D$ of $D\op\tau(\gg,D)$ such that $D\op\tau(\gg,D)=\dot D\op\tau(\gg,D).$ Now for $\dot d\in\dot D,$ there is $d\in D$ and $r\in\tau(\gg,D)$ such that $\dot d=d+r.$ If $x\in\gg,$ we have $[x,\dot d\tilde]=\tau(x,d).$ But  $\dot D$ is a $\gg-$submodule of $\tilde\LL,$ so $[x,\dot d\tilde]=\tau(x,d)\in \dot D\cap\tau(\gg,D)=\{0\}.$  Therefore $\dot D$ is a trivial $\gg-$submodule of $\tilde\LL$ and so $D\op\tau(\gg,D)=\dot D\op\tau(\gg,D)$ is a trivial $\gg-$module. In particular $D$ is a trivial $\gg-$submodule of $\tilde\LL$ and so $\tau(\gg,D)=\{0\}.$ \qed

\begin{Lemma} \label{first} Suppose that $R$ is an irreducible locally finite root system and   $\LL=\op_{\a\in R}\LL_\a$ is an $R-$graded Lie algebra with grading pair $(\gg,\hh).$ Suppose that $\tau:\LL\times\LL\longrightarrow C$ is a $2-$cocycle satisfying $\tau(\LL,\gg)=\{0\}.$ Consider the corresponding central extension $\pi:(\tilde\LL,[\cdot,\cdot\tilde])\longrightarrow \LL$ and suppose  $\tilde\LL$ is  perfect, then $\tilde\LL=\op_{\a\in R}\tilde\LL_\a$ with \begin{equation}\label{decom}\tilde\LL_\a:=\left\{\begin{array}{ll}
\LL_\a&\hbox{if $\a\in R\setminus\{0\}$}\\
\LL_0\op C&\hbox{if $\a=0$}
\end{array}\right.\end{equation} is an $R-$graded Lie algebra with grading pair $(\gg,\hh).$ Moreover, if $R$ is a finite root system, then  the coordinate quadruple of $\LL$ coincides with the coordinate quadruple of $\tilde\LL.$
%
\end{Lemma}
\pf We know  $\tilde\LL=\LL\op ker(\pi)$ and that the corresponding $2-$cocycle $\tau$ satisfies $\tau(\LL,\gg)=\{0\}.$ Since $\tau(\LL,\gg)=\{0\},$ we get that $\gg$ is a subalgebra of $\tilde\LL$ and that  (\ref{decom}) defines a weight space decomposition for $\tilde\LL$ with respect to $\hh.$ So to complete the proof, it is enough to show that $\tilde\LL_0=\sum_{\a\in R\setminus\{0\}}[\tilde\LL_\a,\tilde\LL_{-\a}\tilde].$ For this, we note that  $[\tilde\LL_\a,\tilde\LL_\b\tilde ]\sub\tilde\LL_{\a+\b}$ for $\a,\b\in R$ and so
\begin{eqnarray*}\tilde\LL_0\sub\tilde\LL=[\tilde\LL,\tilde\LL\tilde]
&=&\sum_{\a,\b\in R}[\tilde\LL_\a,\tilde\LL_\b\tilde]\\&=&\sum_{\a,\b;\a+\b=0}[\tilde\LL_\a,\tilde\LL_\b\tilde]+\sum_{\a,\b;\a+\b\neq 0}[\tilde\LL_\a,\tilde\LL_\b\tilde]\\&\sub&\sum_{\a,\b;\a+\b=0}[\tilde\LL_\a,\tilde\LL_\b\tilde]+\sum_{\a,\b;\a+\b\neq 0}\tilde\LL_{\a+\b}.\end{eqnarray*} Now as $\sum_{\a\in R}\tilde\LL_\a$ is direct, we get that $$\displaystyle{\tilde\LL_0=\sum_{\a,\b;\a+\b=0}[\tilde\LL_\a,\tilde\LL_\b\tilde]=[\tilde\LL_0,\tilde\LL_0\tilde]+\sum_{\a\in R\setminus\{0\}}[\tilde\LL_{\a},\tilde\LL_{-\a}\tilde]}.$$
But $\LL_0=\sum_{\a\in R\setminus\{0\}}[\LL_\a,\LL_{-\a}],$ so
\begin{eqnarray*}
[\tilde\LL_0,\tilde\LL_0\tilde]=[\LL_0,\tilde\LL_0\tilde]&=&[\sum_{\a\in R\setminus\{0\}}[\LL_\a,\LL_{-\a}],\tilde\LL_0\tilde]\\
&=&[\sum_{\a\in R\setminus\{0\}}[\LL_\a,\LL_{-\a}\tilde],\tilde\LL_0\tilde]\\
&=&[\sum_{\a\in R\setminus\{0\}}[\tilde\LL_\a,\tilde\LL_{-\a}\tilde],\tilde\LL_0\tilde]\\
&\sub&\sum_{\a\in R\setminus\{0\}}([\tilde\LL_\a,[\tilde\LL_{-\a},\tilde\LL_0\tilde]\tilde]+[\tilde\LL_{-\a},[\tilde\LL_{\a},\tilde\LL_0\tilde]\tilde])\\
&\sub&\sum_{\a\in R\setminus\{0\}}[\tilde\LL_\a,\tilde\LL_{-\a}\tilde].
\end{eqnarray*} So $\tilde\LL_0=\sum_{\a\in R\setminus\{0\}}[\tilde\LL_\a,\tilde\LL_{-\a}].$
This shows that $\tilde\LL$ is an $R-$graded Lie algebra with grading pair $(\gg,\hh).$ Next suppose $R$ is a  finite root system. The Lie algebra epimorphism  $\pi:\tilde\LL\longrightarrow \LL$ induces a Lie algebra epimorphism  $\varphi:\tilde\LL/Z(\tilde\LL)\longrightarrow \LL/Z(\LL)$ mapping $\tilde x+Z(\tilde\LL)$ to $\pi(\tilde x)+Z(\LL)$ for $\tilde x\in\tilde\LL.$ We claim that $\varphi$ is a Lie algebra isomorphism. Suppose that $\tilde x\in\tilde\LL$ and $\pi(\tilde x)\in Z(\LL),$ then for each $\tilde y\in\tilde\LL,$ $\pi([\tilde x,\tilde y\tilde])=[\pi(\tilde x),\pi(\tilde y)]=0$ which implies that $[\tilde x,\tilde y\tilde]\in ker(\pi)\sub Z(\tilde\LL).$ Now it follows that  for each $\tilde y,\tilde z\in \tilde\LL,$ $[\tilde x,[\tilde y,\tilde z\tilde]\tilde]=0,$ and so as $\tilde\LL$ is perfect, we get that $\tilde x\in Z(\tilde\LL),$ therefore $\varphi$ is injective. Now as $\LL$ and $\tilde\LL$ are perfect and $\varphi$ is an isomorphism, we get that  $\LL,$ $\LL/Z(\LL),$ $\tilde\LL$ and $\tilde\LL/Z(\tilde\LL)$ have the same universal central extension, say   $\mathfrak{A}.$ Therefore $\LL$ as well as $\tilde\LL$ are  quotient algebras of  $\mathfrak{A}$ by  subspaces of the center of $\mathfrak{A}.$ Now we are done using \cite[Thm. 420]{ABG1}, \cite[Thm. 5.34]{ABG2}.
\qed

\medskip

Suppose that $\mathfrak{q}=(\fa,*,\cc,f)$ is a coordinate quadruple of type $BC$ and $R$ is an irreducible  locally finite root system of type $BC_I$  for an infinite index set $I.$  Take $\fb:=\fb(\mathfrak{q})$ to be the algebra corresponding to $\mathfrak{q}$ and suppose $\kk$ is a subspace of ${\rm HF}(\fb)$ satisfying the universal property on $\fb.$ Fix a finite subset $I_0$ of $I$ of cardinality greater than $3$ and suppose $\{I_\lam\mid \lam\in \Lam\},$ where  $\Lam$ is an index set containing zero, is the class of all finite subsets of $I$ containing $I_0.$ For $\lam\in\Lam,$ suppose $R_\lam$ is the finite subsystem of $R$ of type $BC_{I_\lam}.$ Next suppose  $\gg=\sum_{\a\in R_{sdiv}}\gg_\a$ is a locally finite split simple Lie algebra of type $C_I$ with splitting Cartan subalgebra $\hh$ and for $\lam\in\Lam,$ take $\gg^\lam:=\sum_{\a\in (R_\lam)_{sdiv}^\times}\gg_\a\op\sum_{\a\in (R_\lam)_{sdiv}^\times}[\gg_\a,\gg_{-\a}].$ One knows that $\gg^\lam$  is a finite dimensional split simple Lie subalgebra of $\gg$ of type $(R_\lam)_{sdiv}$ and $\hh_\lam:=\gg^\lam\cap\hh$ is a splitting  Cartan subalgebra of $\gg^\lam.$ We also recall that $\gg$ is the direct union of $\{\gg^\lam\mid \lam\in\Lam\}.$ Consider the $R-$graded Lie algebra
\begin{equation}\label{decom1}
\LL:=\LL(\mathfrak{q},\kk)=(\gg\ot \aa) \op (\ss\ot \bb)\op (\v\ot\cc)\op \la\fb,\fb\ra
\end{equation} as in Theorem \ref{typebc} in which $\la\fb,\fb\ra:=\{\fb,\fb\}/\kk,$ $\ss$ is an irreducible $\gg-$module equipped with a weight space decomposition with respect to $\hh$ whose set of weights is $R_{lg}\cup\{0\}$  and $\v$ is an irreducible $\gg-$module equipped with a weight space decomposition with respect to $\hh$ whose set of weights is $R_{sh}.$ We also recall that $\ss$ as a vector space is the direct union of the  class $\{\ss^\lam\mid \lam \in \Lam\}$ where  for $\lam\in\Lam,$ $\ss^\lam$ is the irreducible finite  dimensional $\gg^\lam-$module whose set of weights, with respect to $\hh_\lam,$ is $(R_\lam)_{lg}\cup\{0\}$ and $(\ss^\lam)_\a=\ss_\a$ for $\a\in (R^\lam)_{lg},$ also $\v$ as a vector space is the direct union of the  class $\{\v^\lam\mid \lam \in \Lam\}$ where  for $\lam\in\Lam,$ $\v^\lam$ is the  irreducible finite  dimensional $\gg^\lam-$module whose set of weights, with respect to $\hh_\lam,$ is $(R_\lam)_{sh}$ and $(\v^\lam)_\a=\v_\a$ for $\a\in (R^\lam)_{sh}$ (see (\ref{simple-c})).
We next recall from \cite{Y} that for $\lam\in\Lam,$ there is a subalgebra $\dd_\lam$ of $\LL$ with $\dd_0=\la\fb,\fb\ra$ such that  \begin{equation}\label{help5}\begin{array}{l}[\gg^\lam,\dd_\lam]=\{0\},\\\\
\dd_\lam\op (\ss^\lam\dot\ot \bb)=\dd_0\op (\ss^\lam\dot\ot \bb),
\\\\\begin{array}{rl}\LL^\lam:=&(\gg^\lam\dot \ot \aa)\op(\ss^\lam\dot\ot \bb)\op(\v^\lam\dot\ot \cc)\op\dd_\lam\\
=&(\gg^\lam\dot\ot \aa)\op(\ss^\lam\dot\ot \bb)\op(\v^\lam\dot\ot \cc)\op\dd_0\end{array}\end{array}
\end{equation} and that $\LL^\lam$ is a Lie algebra graded by $R_\lam.$
Suppose that  $\pi:(\tilde\LL,[\cdot,\cdot\tilde])\longrightarrow (\LL,[\cdot,\cdot])$ is  a central extension of $\LL.$ As before, we may assume $\tilde\LL=\LL\op ker(\pi),$ $\pi$ is the projection map on $\LL$ and there is a 2-cocycle  $\tau:\LL\times \LL\longrightarrow ker(\pi)$ such that
$$[x_1+z_1,x_2+z_2\tilde]=[x_1,x_2]+\tau(x_1,x_2);\;\;x_1,x_2\in\LL,\;\;z_1,z_2\in ker(\pi).$$ We note that  \begin{eqnarray}\label{help5-1}
{\tilde\LL}^\lam&:=&(\gg^\lam\dot \ot \aa)\op(\ss^\lam\dot\ot \bb)\op(\v^\lam\dot\ot \cc)\op\dd_\lam\op ker(\pi)\nonumber\\
&=&(\gg^\lam\dot\ot \aa)\op(\ss^\lam\dot\ot \bb)\op(\v^\lam\dot\ot \cc)\op\dd_0\op ker(\pi),
\end{eqnarray}
is a central extension of $\LL^\lam.$
Now consider the map  $$\begin{array}{c}\cdot:\gg\times\tilde\LL\longrightarrow \tilde\LL\\x\cdot y\mapsto [x,y\tilde];\;\;x\in\gg,\;y\in\tilde\LL,\end{array}$$
which  defines a $\gg-$module action on  $\tilde\LL.$ One can see that   $\pi$ is a $\gg-$module homomorphism.

Now for each $\lam\in \Lam,$ take
\begin{equation}\label{yahoo3}\mathcal{E}_\lam:=\gg^\lam\op\tau(\gg,\gg).\end{equation}
One can see that  $\mathcal{E}_\lam$ is a $\gg^\lam-$module via the action $``\cdot"$ restricted to $\gg^\lam\times\mathcal{E}_\lam.$ Now as each  finite subset $\{x_1+r_1,\ldots,x_n+r_n\mid x_i\in\gg^\lam,\;r_i\in\tau(\gg,\gg);\;1\leq i\leq n\}$ of $\mathcal{E}_\lam$ is contained in $\gg^\lam\op (\tau(\gg^\lam,\gg^\lam)+\hbox{span}\{r_1,\ldots,r_n\})$ which is a finite dimensional $\gg^\lam-$submodule of the $\gg^\lam-$module $\mathcal{E}_\lam,$ we get that  $\mathcal{E}_\lam$ is a locally finite $\gg^\lam-$module. Therefore   it is completely reducible as $\gg^\lam$ is a finite dimensional split simple Lie algebra. Now as $\tau(\gg,\gg)$ is a $\gg^0-$submodule of $\mathcal{E}_0,$ there is a $\gg^0-$submodule $\dot\gg^0$ of $\mathcal{E}_0$ such that $\mathcal{E}_0=\dot\gg^0\op\tau(\gg,\gg).$
Next we note that $\mathcal{E}_0\sub\mathcal{E}_\lam$ and  for $\lam\in\Lam,$  define
\begin{equation}\label{yahoo4}
\parbox{3.5in}{$\dot\gg^\lam:=$ the $\gg^\lam-$submodule of $\mathcal{E}_\lam$ generated by $\dot\gg^0.$ }
\end{equation}

\begin{Lemma}\label{yahoo1}
(i) Set $\mathcal{E}:=\gg\op\tau(\gg,\gg),$ then $\mathcal{E}$ is both a Lie subalgebra and a $\gg-$submodule  of $\tilde\LL.$ Also the restriction of $\pi$ to $\mathcal{E}$ is both a Lie algebra homomorphism and a $\gg-$module homomorphism.

(ii) $\dot\gg^0$ is a Lie subalgebra of $\tilde\LL$ and the restriction of $\pi$ to $\dot\gg^0$  is both a Lie algebra isomorphism and a $\gg^0-$module isomorphism from $\dot\gg^0$ onto $\gg^0.$ In particular,
$\dot\gg^0$ is a Lie subalgebra of $\tilde\LL$ isomorphic to $\gg^0$ as well as an irreducible $\gg^0-$submodule of $\mathcal{E}_0.$
\end{Lemma}

\pf $(i)$ It is trivial.

$(ii)$ Suppose that $a,b\in\dot\gg^0,$ then since $\gg^0\op\tau(\gg,\gg)=\dot\gg^0\op\tau(\gg,\gg),$ there are unique $x,y\in\gg^0,$ and $r,s\in\tau(\gg,\gg)$ such that $a=x+r$ and  $b=y+s.$ Now as $\dot\gg^0$ is a $\gg^0-$submodule of $\mathcal{E}_0,$ we get that $[a,b\tilde]=[x,b\tilde]\in \dot\gg^0.$ So $\dot\gg^0$ is a Lie subalgebra of $\tilde \LL.$  Next we show that  $\pi_0:=\pi|_{_{\dot\gg^0}}$ is one to one. Suppose that $a,b\in\dot\gg^0$ and $\pi_0(a)=\pi_0(b).$ Since $\dot\gg^0\op\tau(\gg,\gg)=\gg^0\op\tau(\gg,\gg),$ we get that $\pi_0(a)=\pi_0(b)\in \gg^0$ and  that there are unique  $r,s\in\tau(\gg,\gg)$ such that  $a=\pi_0(a)+r$ and $b=\pi_0(b)+s.$  Now as $\pi_0(a)=\pi_0(b)$ and  $\dot\gg^0\cap\tau(\gg,\gg)=\{0\},$ we get that $a-b=r-s=0.$  Now we are done using the fact that $\dot\gg^0\op\tau(\gg,\gg)=\gg^0\op\tau(\gg,\gg).$\qed

%

\medskip

%

\begin{Lemma}\label{yahoo2} Recall (\ref{yahoo4}), we have for $\lam\in \Lam$ that  $\dot\gg^\lam$ is a Lie subalgebra of $\tilde\LL$ and the restriction of $\pi$ to $\dot\gg^\lam$  is both a Lie algebra isomorphism and a $\gg^\lam-$module isomorphism from $\dot\gg^\lam$ to $\gg^\lam.$  In particular, $\dot\gg^\lam$ is a Lie subalgebra of $\tilde\LL$ isomorphic to $\gg^\lam$ and it is an irreducible $\gg^\lam-$submodule of $\tilde\LL$ isomorphic to $\gg^\lam.$ Moreover $\gg^\lam\op\tau(\gg,\gg)=\dot\gg^\lam\op\tau(\gg,\gg).$
\end{Lemma}

\pf We know that $\mathcal{E}_\lam=\gg^\lam\op\tau(\gg,\gg)$ is a locally finite $\gg^\lam-$submodule of $\tilde\LL$ under the action $``\cdot"$ restricted to $\gg^\lam\times \mathcal{E}_\lam$ and that $\tau(\gg,\gg)$ is a submodule of $\gg^\lam\op\tau(\gg,\gg).$ Therefore there is a $\gg^\lam-$submodule $\mathcal{P}$ of $\gg^\lam\op\tau(\gg,\gg)$ such that $\mathcal{E}_\lam=\gg^\lam\op\tau(\gg,\gg)=\mathcal{P}\op\tau(\gg,\gg).$ Then setting $\theta:=\pi|_{_{\mathcal{P}}},$ we get using the same argument as in Lemma \ref{yahoo1}, that  $\theta:\mathcal{P}\longrightarrow \gg^\lam$ is a Lie algebra isomorphism and also a $\gg^\lam-$module isomorphism. Thus as $\gg^\lam$ is equipped with a weight space decomposition $\gg^\lam=(\gg^\lam)_0\op \sum_{\a\in (R_\lam)_{sdiv}^\times}(\gg^\lam)_\a$ with respect to $\hh_\lam,$  we have the weight space decomposition  $\mathcal{P}=\mathcal{P}_0\op\sum_{\a\in (R_\lam)_{sdiv}^\times}{\mathcal P}_\a$ for $\mathcal{P}$ with respect to $\hh_\lam,$ where  $\mathcal{P}_\a:=\theta^{-1}((\gg^\lam)_\a)$ for $\a\in (R_\lam)_{sdiv}.$ This in turn implies that  $(\mathcal{P}_0\op\tau(\gg,\gg))\op\sum_{\a\in (R_\lam)_{sdiv}^\times}{\mathcal P}_\a$ is a weight space decomposition for  $\mathcal{E}_\lam$ with respect to $\hh_\lam.$ We next note $\gg^0$ has a weight space decomposition $\gg^0 =\sum_{\a\in (R_0)_{sdiv}}(\gg^0)_\a$ with respect to $\hh_0$ where  $(\gg^0)_0=\sum_{\a\in (R_0)_{sdiv}^\times}[(\gg^\lam)_\a,(\gg^\lam)_{-\a}]$ and for $\a\in (R_0)_{sdiv}^\times,$ $(\gg^0)_\a=(\gg^\lam)_\a.$ Setting  $\mathcal{Q}:=\theta^{-1}(\gg^0)$ and  $\mathcal{Q}_\a:=\theta^{-1}((\gg^0)_\a)=\mathcal{P}_\a$  for $\a\in (R_0)_{sdiv}\setminus\{0\},$  one gets that $\mathcal{Q}$ is a $\gg^0-$submodule  of $\mathcal{P}$ isomorphic to $\gg^0$ and equipped with the weight space decomposition $\mathcal{Q}=\sum_{\a\in (R_0)_{sdiv}^\times}\mathcal{Q}_\a\op\sum_{\a\in (R_0)_{sdiv}^\times}[\mathcal{Q}_\a,\mathcal{Q}_{-\a}]$  with respect to $\hh_0.$ Also   $\sum_{\a\in (R_0)_{sdiv}^\times}\mathcal{Q}_\a\op(\sum_{\a\in (R_0)_{sdiv}^\times}[\mathcal{Q}_\a,\mathcal{Q}_{-\a}]\op \tau(\gg,\gg))$  and $\mathcal{E}_0=\mathcal{Q}\op \tau(\gg,\gg)$  is a weight space decomposition of $\mathcal{E}_0$ with respect to $\hh_0.$  Now $\dot\gg^0$ is a nontrivial finite dimensional irreducible  $\gg^0-$submodule of $\mathcal{E}_0$ isomorphic to $\gg^0$ and so by [Y,Theorem], $\dot\gg^\lam$ is a $\gg^\lam-$submodule of $\mathcal{E}_\lam$ isomorphic to $\gg^\lam.$
On the other hand we know that $\theta:\mathcal{E}_\lam\longrightarrow \gg^\lam$ is a $\gg^\lam-$module homomorphism. Now as $\gg^\lam$ and $\dot\gg^\lam$ are irreducible $\gg^\lam-$modules and $\theta(\dot\gg^0)=\pi(\dot\gg^0)\neq0,$ one gets that the restriction of  $\pi$ to $\dot\gg^\lam$ is a $\gg^\lam-$module isomorphism from $\dot\gg^\lam$ onto $\gg^\lam$ which in turn implies that $\dot\gg^\lam\op\tau(\gg,\gg)=\gg^\lam\op\tau(\gg,\gg)$ and  that $\pi\mid_{_{\dot\gg_\lam}}:\dot\gg^\lam\longrightarrow \gg^\lam$ is also a  Lie algebra isomorphism. \qed

\begin{cor}
For $\lam,\mu\in\Lam$ with $\lam\prec \mu,$ we have   $\dot\gg^\lam\sub\dot\gg^\mu,$ in particular $\cup_{\lam\in\Lam}\dot\gg^\lam$ is a subalgebra of $\tilde\LL$ and also a $\gg-$submodule of $\tilde\LL.$ Also setting  $\dot\gg$ to be the direct union of $\{\dot\gg^\lam\mid \lam\in\Lam\}, $ $\pi\mid_{\dot\gg} $ is both a Lie algebra isomorphism and a  $\gg-$module isomorphism  from $\dot\gg$ to $\gg.$ Moreover,  we have $\gg\op\tau(\gg,\gg)=\dot\gg\op\tau(\gg,\gg).$
\end{cor}

\pf Since $\dot\gg^0\sub\dot\gg^\lam\cap\dot\gg^\mu$ and $\gg^\lam\sub\gg^\mu,$ we get that $\dot\gg^\lam\sub\dot\gg^\mu.$ Now using  Lemmas \ref{yahoo1} and  \ref{yahoo2}, we are done.\qed

\bigskip

Recall (\ref{decom1}) and suppose $\ii$ is an index set containing  zero and fix a basis $\{a_i\mid i\in \ii \}$ with $a_0=1$ for $\aa,$ also fix a basis $\{b_j\mid j\in \jj\}$ for $\bb$ and a basis $\{c_t\mid t\in \T\}$ for $\cc.$ For $i\in\ii,$ $j\in\jj$ and $t\in\T,$ set \begin{equation}\label{new1}\gg_i:=\gg\ot a_i,\;\;\ss_j:=\ss\ot b_j,\;\;\v_t:=\v\ot c_t\end{equation} which are $\gg-$submodules of $\LL;$ also for $\lam\in\Lam,$ set
\begin{equation}\label{new2}\gg^\lam_i:=\gg^\lam\ot a_i,\;\;\ss^\lam_j:=\ss^\lam\ot b_j,\;\;\v^\lam_t:=\v^\lam\ot c_t.\end{equation}
Now suppose $\m$ is one of the  $\gg-$submodules of $\LL$ in the class $\{\gg_i,\ss_j,\v_t\mid i\in\ii\setminus\{0\},\;j\in\jj,\;t\in\T\}$ and consider  $\tilde\m:=\m\op\tau(\gg,\m).$ Then $\tilde\m$ is a $\gg-$submodule of $\tilde\LL.$ We know that the vector space $\m$ is the direct union of a class $\{\m^\lam\mid \lam\in \Lam\}$ in which  each $\m^\lam$ is a finite dimensional irreducible  $\gg^\lam-$submodule of $\LL^\lam$ equipped with the weight spaced decomposition $\m^\lam=\op_{\gamma\in \Gamma_\lam }(\m^\lam)_\gamma$ where 
\begin{equation}\label{new7}
\parbox{5in}{\begin{itemize}
\item $\Gamma_\lam\sub R_\lam,$
\item For $*\in\{sh,lg,ex\},$ $\Gamma_\lam\setminus\{0\}=(R_\lam)_*$ if and only if  $\Gamma_0\setminus\{0\}=(R_0)_*,$
\item $\Gamma_\lam\sub \Gamma_\mu;\;\;\; \lam\prec\mu,$
\item $(\m^\lam)_\gamma=(\m^\mu)_\gamma;$ $\lam\prec\mu, $ $\gamma\in \Gamma_\lam\setminus\{0\}.$
\end{itemize}}
\end{equation}  For  $\lam\in \Lam,$ set $\tilde\m^\lam:=\m^\lam\op\tau(\gg,\m),$ then   $\tilde\m^\lam$ is a $\gg^\lam-$submodule of $\tilde\LL$ under the action $``\cdot"$ restricted to $\gg^\lam\times\tilde\LL.$
Now if  $\{m_i+r_i\mid 1\leq i\leq n,\; m_i\in\m^\lam,r_i\in \tau(\m,\gg)\}$ is a finite subset of $\tilde\m^\lam,$ we see that  $\{m_i+r_i\mid 1\leq i\leq n\}\sub\m^\lam+\tau(\gg^\lam,\m^\lam)+\hbox{span}_\bbbf\{r_1,\ldots,r_n\}$ which is a finite dimensional $\gg^\lam-$submodule of $\tilde\LL.$ This  means that $\tilde\m^\lam$ is a locally finite $\gg^\lam-$module and so it is completely reducible as $\gg^\lam$ is finite dimensional simple Lie algebra.

Next we note that  $\tau(\gg,\m)$ is a $\gg^0-$submodule of the locally finite $\gg^0-$module $\tilde\m^0$, so  there is a $\gg^0-$submodule $\dot\m^0$ of $\tilde\m^0$ such that $\tilde\m^0=\dot\m^0\op\tau(\gg,\m).$

Set
\begin{equation}\label{rest1}
\dot\m^\lam:=\hbox{$\gg^\lam-$submodule of $\tilde\m^\lam$ generated by $\dot\m^0;$}\;\; \lam\in \Lam.\end{equation}

\begin{Lemma}\label{uce1}
(i) For $\lam\in\Lam,$ the restriction of $\pi$ to $\dot\m^\lam$ is a $\gg^\lam-$module isomorphism   from $\dot\m^\lam$ onto $\m^\lam$   and $\tilde\m^\lam=\dot\m^\lam\op\tau(\gg,\m).$

(ii) For $\lam\prec\mu,$ we have $\dot\m^\lam\sub\dot\m^\mu,$ in particular $\dot\m,$ the direct union of $\{\dot\m^\lam\mid \lam\in\Lam\},$  is a $\gg-$submodule of $\tilde\LL.$ Also   the restriction of $\pi$ to $\dot\m$ is a $\gg-$module isomorphism from $\dot\m$ onto $\m$ and $\tilde\m=\dot\m\op\tau(\gg,\m)=\m\op\tau(\gg,\m).$

(iii) If $x\in\dot\m$ and $\pi(x)\in\m^\lam$ for some $\lam\in\Lam,$ then $x\in\dot\m^\lam.$
\end{Lemma}
\pf $(i)$ We first note that since $\pi$  is  a $\gg-$module homomorphism, the restriction of $\pi$ to $\tilde\m$ is a $\gg-$module homomorphism. Now as  $\m^0\op\tau(\gg,\m)=\dot\m^0\op\tau(\gg,\m)$ and that $\m^0$ is an irreducible  $\gg^0-$module, it is immediate that the restriction of $\pi$ to $\dot\m^0$ is a $\gg^0-$module isomorphism  from $\dot\m^0$ onto $\m^0.$
Now suppose that $0\prec\lam,$ since $\tilde\m^\lam$ is a completely reducible $\gg^\lam-$module and $\tau(\gg,\m)$ is a $\gg^\lam-$submodule of $\tilde\m^\lam,$ one finds a $\gg^\lam-$submodule of $\mathcal{N}$ of $\tilde\m^\lam$ such that $\tilde\m^\lam=\mathcal{N}\op\tau(\gg,\m).$ Therefore $\theta:=\pi|_{_{\mathcal{N}}}:\mathcal{N}\longrightarrow \m^\lam$ is a $\gg^\lam-$module isomorphism. We know that $\m^\lam$ has a weight space decomposition $\m^\lam=\op_{\a\in\Gamma_\lam}(\m^\lam)_\a$  with respect to $\hh_\lam$ and that  $\m^0$ has a weight space decomposition $\m^0=\op_{\a\in\Gamma_0}(\m^0)_\a $ with respect to $\hh_0$ such that
\begin{equation}\label{help3}
\begin{array}{l}
\Gamma_0\sub R_0,\;\;\;\;\Gamma_\lam\sub R_\lam,\\
\hbox{For $*\in\{sh,lg,ex\},$ $\Gamma_\lam\setminus\{0\}=(R_\lam)_*$ if and only if  $\Gamma_0\setminus\{0\}=(R_0)_*,$}\\
(\m^0)_\a=(\m^\lam)_\a \hbox{ for } \a\in \Gamma_0\setminus\{0\},
\end{array}
\end{equation}
(see  (\ref{new7})). Now since  $\m^0$ is a $\gg^0-$submodule of $\m^\lam$ and $\theta $ is a $\gg^\lam-$module isomorphism, $\mathcal{N}^0:=\theta^{-1}(\m^0)\sub\m^0\op\tau(\gg,\m)=\tilde\m^0$ is a $\gg^0-$submodule  of $\tilde\m^0.$ Also $\mathcal{N}$ is a $\gg^\lam-$module equipped with the weight space decomposition  $\mathcal{N}=\op_{\a\in\Gamma_\lam}\mathcal{N}_\a$ with respect to $\hh_\lam,$ where for $\a\in \Gamma_\lam,$ $\mathcal{N}_\a:=\theta^{-1}((\m^\lam)_\a),$  and that  $\mathcal{N}^0$ has a weight space decomposition $\mathcal{N}^0=\op_{\a\in\Gamma_0}(\mathcal{N}^0)_\a,$  with respect to $\hh_0,$ where for $\a\in \Gamma_0,$ $(\mathcal{N}^0)_\a:=\theta^{-1}((\mathcal{M}^0)_\a).$ Therefore $\tilde\m^\lam$  has a weight space decomposition $\tilde\m^\lam=\op_{\a\in\Gamma_\lam\cup\{0\}}(\tilde\m^\lam)_\a,$ with respect to $\hh_\lam$ where
\begin{equation}\label{help1}(\tilde\m^\lam)_\a=\left\{\begin{array}{ll}(\m^\lam)_\a& \hbox{if $\a\in\Gamma_\lam\setminus\{0\},$}\\
(\m^\lam)_0+\tau(\gg,\m)&\hbox{if $\a=0$ and  $0\in\Gamma_\lam,$ }\\
\tau(\gg,\m)&\hbox{if $\a=0$ and  $0\not\in\Gamma_\lam.$ }
\end{array}\right.\end{equation}
Also $\tilde\m^0=\mathcal{N}^0\op\tau(\gg,\m)$
and $\tilde\m^0$ is equipped with the weight space decomposition
$\tilde\m^0=\op_{\a\in\Gamma_0\cup\{0\}}(\tilde\m^0)_\a$ with respect to $\hh_0,$ where \begin{equation}\label{help2}(\tilde\m^0)_\a=\left\{\begin{array}{ll}(\m^0)_\a& \hbox{if $\a\in\Gamma_0\setminus\{0\},$}\\
(\m^0)_0+\tau(\gg,\m)&\hbox{if $\a=0$ and  $0\in\Gamma_0,$ }\\
\tau(\gg,\m)&\hbox{if $\a=0$ and  $0\not\in\Gamma_0.$ }
\end{array}\right.\end{equation}
Now (\ref{help3})-(\ref{help2}) together with [Y, Theorem] imply that the $\gg^\lam-$submodule $\dot\m^\lam$ of $\tilde\m^\lam$ generated by $\dot\m^0$ is a $\gg^\lam-$submodule of $\tilde\m^\lam$ isomorphic to $\m^\lam. $  This together with the facts that $\pi(\dot\m^\lam)\sub\m^\lam,$
$\pi(\dot\m^0)\neq\{0\},$
and $\dot\m^\lam$ as well as $\m^\lam$ are   irreducible $\gg^\lam-$modules, implies that the restriction of $\pi$ to $\dot\m^\lam$ is a $\gg^\lam-$module isomorphism from $\dot\m^\lam$ to $\m^\lam.$ In particular, we get that  $\tilde\m^\lam=\m^\lam\op\tau(\gg,\m)=\dot\m^\lam\op\tau(\gg,\m).$

$(ii)$ This is easy to see using Part $(i).$

$(iii)$ Take $\mu\in\Lam$ to be such that $x\in\dot\m^\mu.$ We know that there is  $\nu\in\Lam$ with $\lam\prec\nu$ and $\mu\prec\nu.$ Since the restriction of $\pi$ to $\dot\m^\lam$ is a  $\gg^\lam-$module isomorphism from $\dot\m^\lam$ onto $\m^\lam,$ one finds $y\in\dot\m^\lam$ such that $\pi(x)=\pi(y).$ So $x,y\in\dot\m^\nu$ and $\pi(x)=\pi(y).$ But the restriction of $\pi$ to $\dot\m^\nu$ is a  $\gg^\nu-$module isomorphism from $\dot\m^\nu$ onto $\m^\nu,$ therefore $x=y\in\dot\m^\lam.$
\qed

\bigskip
Consider  (\ref{new1}) and (\ref{new2}) and  identify $\gg\ot 1$ with $\gg.$ Using  Lemmas \ref{uce1}, \ref{yahoo1} and \ref{yahoo2},  if  $i\in\i\setminus,$ $j\in\jj$ and $t\in\T,$ for $\lam\in\Lam,$ one finds irreducible $\gg^\lam-$submodules $\dot\gg_i^\lam,$ $\dot\ss_j^\lam$ and $\dot\v_t^\lam$ of $\tilde\LL$ such that
 $\dot\gg_i^\lam$ is isomorphic to $\gg_i^\lam,$ $\dot\ss_j^\lam$ is isomorphic to $\ss_j^\lam$ and $\dot\v_t^\lam$ is isomorphic to $\v_t^\lam.$
Moreover
\begin{equation}\label{new4}
\parbox{4in}{
\begin{itemize}
\item $\dot\gg_i^\lam$ is the $\gg^\lam-$submodule of $\tilde\LL$ generated by $\dot\gg_i^0,$
\item  $\dot\ss_j^\lam$ is the $\gg^\lam-$submodule of $\tilde\LL$ generated by $\dot\ss_j^0,$
\item $\dot\v_t^\lam$ is the $\gg^\lam-$submodule of $\tilde\LL$ generated by $\dot\v_t^0.$
\end{itemize}}\end{equation}

Also setting $\dot\gg_i:=\underrightarrow{\lim}_{\lam\in\Lam}\dot\gg_i^\lam,$  $\dot\ss_j:=\underrightarrow{\lim}_{\lam\in\Lam}\dot\ss_j^\lam$ and
$\dot\v_t:=\underrightarrow{\lim}_{\lam\in\Lam}\dot\v_t^\lam,$
 $\dot\gg_i$ is isomorphic to $\gg_i,$ $\dot\ss_j$ is isomorphic to $\ss_j$ and $\dot\v_t$ is isomorphic to $\v_t.$ Also we have
\begin{equation}\label{new5}
\parbox{4.3in}{
\begin{itemize}
\item $\gg_i\op\tau(\gg,\gg_i)=\dot\gg_i\op\tau(\gg,\gg_i),\;$  $\gg_i^\lam\op\tau(\gg,\gg_i)=\dot\gg^\lam_i\op\tau(\gg,\gg_i),$
\item  $\ss_j\op\tau(\gg,\ss_j)=\dot\ss_j\op\tau(\gg,\ss_j),\;$   $ \ss_j^\lam\op\tau(\gg,\ss_j)=\dot\ss^\lam_j\op\tau(\gg,\ss_j),$
\item $\v_t\op\tau(\gg,\v_t)=\dot\v_t\op\tau(\gg,\v_t),\;$  $ \v_t^\lam\op\tau(\gg,\v_t)=\dot\v_t^\lam\op\tau(\gg,\v_t).$
\end{itemize}}\end{equation}

\begin{Lemma}\label{uce9}
Consider (\ref{help5}), there is a subspace $\dot\dd$ of $[\LL^0,\LL^0\tilde]\cap (\dd_0\op ker(\pi))$ such that $\pi(\dot\dd)\sub\dd_0,$  $\pi|_{\dot\dd}:\dot\dd\longrightarrow \dd_0$ is a linear isomorphism, $[\gg^0,\dot\dd\tilde]=\{0\}$ and for $\lam\in\Lam,$ $[\gg^\lam,\dot\dd\tilde]\sub\dot\sum_{j\in\jj}\dot\ss_j^\lam.$
\end{Lemma}

\pf We note that $\dd_0\sub\LL^0=[\LL^0\,\LL^0],$ so for $d\in\dd_0,$ there  are $n\in\bbbn\setminus\{0\},$ $x_i,y_i\in\LL^0$ such that $d=\sum_{i=1}^n[x_i,y_i].$ So we have $\sum_{i=1}^n[x_i,y_i\tilde]=d+\sum_{i=1}^n\tau(x_i,y_i)\in[\LL^0,\LL^0\tilde]\cap (\dd_0+ker(\pi)).$ Also $\pi(\sum_{i=1}^n[x_i,y_i\tilde])=d.$
Therefore there is a subspace $\dot\dd$ of $[\LL^0,\LL^0\tilde]\cap (\dd_0\op ker(\pi))$ such that $\pi|_{\dot\dd}:\dot\dd\longrightarrow \dd_0$ is a linear isomorphism. Now using Lemma \ref{uce3} together with the fact that $\pi\mid_{_{\LL^0\op ker(\pi)}}:\LL^0\op ker(\pi)\longrightarrow \LL^0$ is  a central extension of $\LL^0,$ we get that $[\gg^0,\dot\dd\tilde]\sub[\gg^0,\dd_0\tilde]=\{0\}.$ Next suppose $\lam\in\Lam,$ $x\in\gg^\lam$ and $\dot d\in\dot\dd\sub\dd_0\op ker(\pi)\sub\dd_\lam+\sum_{j\in \jj}\ss_j^\lam+ker(\pi)$ (see (\ref{help5})). So as $\LL^\lam\op ker(\pi)$ is a central extension for $\LL^\lam,$ using Lemma \ref{uce3} together with (\ref{new5}), we have
\begin{eqnarray*}[x,\dot d\tilde]\in[x,\dd_\lam+\sum_{j\in \jj}\ss_j^\lam+ker(\pi)\tilde]\sub[x,\dd_\lam+\sum_{j\in \jj}\ss_j^\lam\tilde]
&\sub&[x,\sum_{j\in \jj}\ss_j^\lam\tilde]\\
&\sub&[x,\sum_{j\in \jj}\dot\ss_j^\lam\tilde]\\
&\sub&\sum_{j\in \jj}\dot\ss_j^\lam.
\end{eqnarray*}This completes the proof.\qed

\begin{Lemma}\label{uce2}
$\sum_{i\in I}\dot\gg_i+\sum_{j\in J}\dot\ss_j+\sum_{t\in \T}\dot\v_t+ \dot\dd$ is a direct sum.
\end{Lemma}
\pf Suppose that $\sum_i\dot x_i\in \sum_{i\in \ii}\dot\gg_i,\sum_{j}\dot y_j\in \sum_{j\in \j}\dot\ss_j,\sum_{t}\dot z_t\in \sum_{t\in \T}\dot\v_t,\dot d\in \dot\dd$ and $\sum_i\dot x_i+\sum_{j}\dot y_j+\sum_{t}\dot z_t+\dot d=0,$ then  we have
$$0=
\pi(\sum_i\dot x_i+\sum_{j}\dot y_j+\sum_{t}\dot z_t+\dot d)=\sum_i\pi(\dot x_i)+\sum_{j}\pi(\dot y_j)+\sum_{t}\pi(\dot z_t)+\pi(\dot d).$$ Now as  $\pi(\dot\gg_i)=\gg_i,$ $\pi(\dot\ss_j)=\ss_j,$  $\pi(\dot\v_t)=\v_t$ and  $\sum_{i\in \ii}\gg_i+\sum_{j\in \j}\ss_j\op\sum_{t\in \T}\v_t+\dd$ is direct, we get that $\pi(\dot x_i)=0,\pi(\dot y_j)=0,\pi(\dot z_t)=0$ and $\pi(\dot d)=0.$ But for $i\in\ii,$ $j\in \jj$ and $t\in\T,$ $\pi\mid_{_{\dot\gg_i}}, \pi\mid_{_{\dot\ss_j}},$  $\pi\mid_{_{\dot\v_t}}$ and $\pi\mid_{_{\dot\dd}}$ are isomorphism, so $\dot x_i=0,\dot y_j=0,\dot z_t=0,\dot d=0$ ($i\in\ii,$ $j\in \jj$ and $t\in\T$).
\qed
\bigskip

Now set
\begin{equation}\label{help4}\dot\LL:=\bigoplus_{i\in \ii}\dot\gg_i\op\bigoplus_{j\in \jj}\dot\ss_j\op\bigoplus_{t\in \T}\dot\v_t\op\dot\dd.
\end{equation} Using the same argument as in Lemma \ref{uce2}, $\dot\LL\cap ker(\pi)=\{0\}$ and so $\tilde\LL=\dot\LL\op ker(\pi).$

 Now set \begin{equation}\label{new}\pi_1:\tilde\LL\longrightarrow \dot\LL\andd\pi_2:\tilde\LL\longrightarrow ker(\pi)\end{equation} to be the projective  maps on $\dot\LL$ and $ker(\pi)$ respectively and for $x,y\in\dot\LL,$ define
\begin{equation}\label{yahoo8}
[x,y\dot]:=\pi_1([x,y\tilde])\andd \dot\tau(x,y)=\pi_2([x,y\tilde]).\end{equation}
Then we get that $(\dot\LL,[\cdot,\cdot\dot])$ is a Lie algebra and by Lemma \ref{yahoo2}, $\dot\gg^\lam,$ $\lam\in\Lam,$ is a subalgebra of $\dot\LL.$  Also  $\dot\tau:\dot\LL\times\dot\LL\longrightarrow ker(\pi)$ is a 2-cocycle. Moreover consulting Lemma \ref{uce9}, we get that  $\dot\LL$ is a $\gg-$submodule of $\tilde\LL$ which implies that
\begin{equation}\label{uce4}\dot\tau(\dot\gg,\dot\LL)=\{0\}.\end{equation}

For each $x\in\tilde\LL,$ there are unique $\dot\ell_x\in\dot\LL,$ $\ell_x\in\LL$ and $e_x,f_x\in ker(\pi)$ such that $x=\dot\ell_x+e_x=\ell_x+f_x.$  Also we note that
\begin{equation}\label{yahoo6}
\begin{array}{l}
\ell_{\dot\ell_x}=x \andd f_{\dot\ell_x}=-e_x;\;\;\;x\in\LL\\
\dot\ell_{\ell_y}=y\andd e_{\ell_y}=-f_y;\;\;\;y\in\dot\LL,
\end{array}
\end{equation}
(see (\ref{new5})). For $\lam\in\Lam,$ set
$$\dot\LL^\lam:=\sum_{i\in \ii}\dot\gg_i^\lam\op\sum_{j\in \jj}\dot\ss_j^\lam\op\sum_{t\in \T}\dot\v_t^\lam\op\dot\dd.$$ So considering (\ref{help5-1}), we have \begin{equation}\label{help6}
\tilde\LL^\lam:=\LL^\lam\op  ker(\pi)=\dot\LL^\lam\op ker(\pi).\end{equation}
Note that for each $\lam\in\Lam,$ $\tilde\LL^\lam$ is a Lie subalgebra of $\tilde\LL$ and $\tilde\LL$ is the direct union of $\{\tilde\LL^\lam\mid \lam\in\Lam\}.$ In the following, we show that for $\lam\in\Lam,$ $\dot\LL^\lam$ is a Lie subalgebra of $\dot\LL$ and that $\dot\LL$ is the direct union of $\{\dot\LL^\lam\mid \lam\in\Lam\}.$
\begin{Lemma}\label{uce6}

(i)
$\pi|_{\dot\LL}$ is a Lie algebra isomorphism from $(\dot\LL,[\cdot,\cdot\dot])$ to  $(\LL,[\cdot,\cdot]).$ Also for each $\lam\in\Lam,$ $\dot\LL^\lam$ is a Lie subalgebra of $(\dot\LL,[\cdot,\cdot\dot])$ isomorphic to $\LL^\lam$ and $\dot\LL$ is the direct union of $\{\dot\LL^\lam\mid\lam\in\Lam\}.$ 

(ii) Recall (\ref{new}), for $\lam\in\Lam,$ ${\pi_1}|_{_{\tilde\LL^\lam}}:\tilde\LL^\lam\longrightarrow \dot\LL^\lam$ is a central extension of $\dot\LL^\lam$ with corresponding 2-cocycle $\dot\tau\mid_{_{\dot\LL^\lam\times\dot\LL^\lam}}$ satisfying $\dot\tau(\dot\gg^\lam,\dot\LL^\lam)=\{0\}.$

(iii) For $\lam\in\Lam,$ there is a subalgebra $\dot\dd_\lam$ of $\dot\LL^\lam$ with $\pi(\dot\dd_\lam)=\dd_\lam$ such that $\dot\dd_\lam$ is a trivial  $\dot\gg^\lam-$submodule of $\dot\LL^\lam$ and \begin{eqnarray*}\dot\LL^\lam&=&\sum_{i\in \ii}\dot\gg_i^\lam\op\sum_{j\in \jj}\dot\ss_j^\lam\op\sum_{t\in \T}\dot\v_t^\lam\op\dot\dd\\
&=&\sum_{i\in \ii}\dot\gg_i^\lam\op\sum_{j\in \jj}\dot\ss_j^\lam\op\sum_{t\in \T}\dot\v_t^\lam\op\dot\dd_\lam.
\end{eqnarray*}
\end{Lemma}
\pf 
($i$) We fix $x,y\in\dot\LL$ and show that $\pi([x,y\dot])=[\pi(x),\pi(y)]=[\ell_x,\ell_y].$ Using (\ref{yahoo6}), we have
\begin{eqnarray*}\pi([ x, y\dot])=\pi(\pi_1([ x, y\tilde]))&=&\pi(\pi_1([\ell_x,\ell_y]+\tau(\ell_x,\ell_y)))\\&=&\pi(\pi_1(\dot\ell_{[\ell_x,\ell_y]}+e_{[\ell_x,\ell_y]}+\tau(\ell_x,\ell_y)))\\&=&
\pi(\dot\ell_{[\ell_x,\ell_y]})\\&=&\ell_{\dot\ell_{[\ell_x,\ell_y]}}=[\ell_x,\ell_y].\end{eqnarray*} This means that the restriction of  $\pi$ to $\dot\LL$ is a Lie algebra homomorphism. But $\tilde\LL=\LL\op ker(\pi)=\dot\LL\op ker(\pi)$ which in turn implies that $\pi$ restricted to $\dot\LL$ is an isomorphism from $\dot\LL$ onto $\LL.$
Next suppose $\lam\in\Lam$ and  $x,y\in\dot\LL^\lam,$ then $\ell_x,\ell_y\in\LL^\lam.$  Also  $[x,y\dot]\in\dot\LL,$ and $[x,y\dot]=[x,y\tilde]-\dot\tau(x,y)=[\ell_x,\ell_y]+\tau(\ell_x,\ell_y)-\dot\tau(x,y)\in\LL^\lam+ ker(\pi)=\tilde\LL^\lam.$ Therefore we get $[x,y\dot]\in\dot\LL\cap\tilde\LL^\lam=\dot\LL^\lam$ which shows that $\dot\LL^\lam$ is a subalgebra of $\dot\LL.$ Now as the restriction of $\pi$ to $\dot\LL$ is a Lie algebra isomorphism from $\dot\LL$ to $\LL,$ we get using (\ref{help6}) that for $\lam\in\Lam,$ the restriction of $\pi$ to $\dot\LL^\lam$ is a Lie algebra isomorphism from $\dot\LL^\lam$ to $\LL^\lam.$ Now  consider (\ref{help5}) and set $\dot\dd_\lam:=\dot\LL^\lam\cap \pi^{-1}(\dd_\lam).$  We note that $\pi:\tilde\LL^\lam\longrightarrow \LL^\lam$ is a central extension, so using Lemma \ref{uce3}, we have $[\dot\gg^\lam,\dot\dd_\lam\tilde]=[\gg^\lam,\dd_\lam\tilde]=0.$

$(ii)$ Since $\tilde\LL^\lam=\dot\LL^\lam\op ker(\pi)$,  it is immediate that ${\pi_1}|_{_{\tilde\LL^\lam}}$ is a central extension of $\dot\LL^\lam.$  Also we note that for $x\in\dot\gg^\lam$ and $y\in\dot\LL^\lam,$ $[\ell_x,y\tilde]\in\dot\LL^\lam$ as $\dot\LL^\lam$ is a $\gg^\lam-$submodule of $\tilde\LL.$ Thus $\dot\tau(x,y)=\pi_2([x,y\tilde])=\pi_2([\ell_x,y\tilde])=0.$\qed

\bigskip

Now using the same notation as in the text and regard Lemmas \ref{uce9} and \ref{uce6}, we summarize our information as follows:
We have $$\begin{array}{l}
\dot\dd\sub[\LL^0,\LL^0\tilde]=[\dot\LL^0,\dot\LL^0\tilde],\\\\
\tilde\LL=\dot\LL\op ker(\pi)=\sum_{i\in\i}\dot\gg_i\op\sum_{j\in\jj}\dot\ss_j\op\sum_{t\in\T}\dot\v_t\op\dot\dd\op ker(\pi),\\\\
\begin{array}{ll}\tilde\LL^\lam=\dot\LL^\lam\op ker(\pi)&=\sum_{i\in\i}\dot\gg^\lam_i\op\sum_{j\in\jj}\dot\ss^\lam_j\op\sum_{t\in\T}\dot\v^\lam_t\op\dot\dd\op ker(\pi)\\
&=\sum_{i\in\i}\dot\gg^\lam_i\op\sum_{j\in\jj}\dot\ss^\lam_j\op\sum_{t\in\T}\dot\v^\lam_t\op\dot\dd_\lam\op ker(\pi),\end{array}
\end{array}$$ $(\lam\in\Lam).$ Also  $\pi_1:\tilde\LL\longrightarrow \dot\LL$ is a central extension of $\dot\LL$ with corresponding 2-cocycle $\dot\tau$ satisfying $\dot\tau(\dot\gg,\dot\LL)=\{0\}.$ So  without loss of generality, from now on we assume that $\tilde\LL=\LL\op ker(\pi)$ and that $\LL$ is a $\gg-$submodule of $\tilde\LL$ and that \begin{equation}\label{yahoo9}\tau(\gg,\LL)=\{0\} \andd
\dd_0\sub[\LL^0,\LL^0\tilde]=[\tilde\LL^0,\tilde\LL^0\tilde].\end{equation}

\begin{Lemma}\label{uce8}
For $\lam\in\Lam,$ $\LL^\lam\sub [\LL^\lam,\LL^\lam\tilde].$ Also  $\LL\sub [\LL,\LL\tilde].$
\end{Lemma}
\pf We know from (\ref{help5}) that
\begin{equation}\label{decom2}\LL^\lam=(\gg^\lam\dot\ot \aa)\op(\ss^\lam\dot\ot \bb)\op(\v^\lam\dot\ot \cc)\op\dd_\lam=(\gg^\lam\dot\ot \aa)\op(\ss^\lam\dot\ot \bb)\op(\v^\lam\dot\ot \cc)\op\dd_0.\end{equation}
and that  $\LL=(\gg\ot \aa)\op(\ss\ot \bb)\op(\v\ot \cc)\op \dd_0.$  We also know that 
 the restriction $\pi$ to $\LL^\lam\op ker(\pi)$ is a central extension of $\LL^\lam$ with corresponding 2-cocycle $\tau$ satisfying $\tau(\gg^\lam,\LL^\lam)=\{0\}.$ Thus it follows from   \cite[Pro. 5.23]{ABG2} that the  summands $\gg^\lam\ot \aa,$ $\ss^\lam\ot \bb,$ $\v^\lam\ot \cc$ and $\dd^\lam$ are orthogonal with respect to $\tau$ and that  for  $x,y\in\gg^\lam,$ $a\in \aa,$ $s\in \ss,$ $b\in\bb,$ $v\in\v$ and  $c\in\cc,$ we have
$$\begin{array}{l}\;[x\ot 1,y\ot a\tilde]=[x\ot 1,y\ot a]+\tau(x\ot 1,y\ot a)=[x\ot 1,y\ot a]=[x,y]\ot a,\\
\;[x\ot 1,s\ot b\tilde]=[x\ot 1,s\ot b]=[x,s]\ot b,\\
\;[x\ot 1,v\ot c\tilde]=[x\ot 1,v\ot c]=xv\ot c.
\end{array}$$ 
This together with  the fact that $\gg,\ss$ and $\v$ are irreducible finite dimensional $\gg-$modules, implies that $(\gg^\lam\dot\ot \aa)\op(\ss^\lam\dot\ot \bb)\op(\v^\lam\dot\ot \cc)\sub[\tilde\LL^\lam,\tilde\LL^\lam\tilde].$ Now contemplating (\ref{decom2}) and (\ref{yahoo9}), we are done.
\qed

\begin{Theorem}\label{main'} Suppose that $I$ is  an infinite index set, $R$ is an irreducible locally finite root system of type $BC_I$ and $\mathfrak{q}:=(\fa,*,\cc,f)$ is a coordinate quadruple of type $BC.$ Take   $\fb:=\fb(\mathfrak{q})$ and suppose $\kk$ is a subspace of ${\rm HF}(\fb)$ satisfying the  uniform property on $\fb.$ Set $\la\fb,\fb\ra:=\{\fb,\fb\}/\kk$ and consider the $R-$graded Lie algebra   $\LL:=(\gg\ot \aa)\op(\ss\ot \bb)\op(\v\ot \cc)\op\la\fb,\fb\ra.$ Suppose that $\tau:\LL\times \LL\longrightarrow E$ is a 2-cocycle and consider the corresponding  central extension $\tilde\LL:=\LL\op E$ as well as  the canonical projection map $\pi:\tilde\LL\longrightarrow \LL.$ If $\tilde\LL$ is perfect, then $\tilde\LL$ is an $R-$graded Lie algebra with the same coordinate  quadruple $(\fa,*,\cc,f).$ Also there is a subspace $\kk_0$ of ${\rm HF}(\fb)$ satisfying the uniform property on $\fb$ such that $\tilde\LL$ can be identified with $(\gg\ot \aa)\op(\ss\ot \bb)\op(\v\ot \cc)\op(\{\fb,\fb\}/\kk_0)$  where setting $\la b,b'\ra_c:=\{b,b'\}+\kk_0,$ the  Lie bracket on $\tilde\LL$ is given by \begin{equation}\label{probc-fin}
\begin{array}{l}
\;[x\ot a,y\ot a'\tilde]=[x,y]\ot\frac{1}{2}(a\circ a')+ (x\circ y)\ot\frac{1}{2}[a,a']+tr(xy)\la a,a'\ra_c,\vspace{1mm}\\
\;[x\ot a,s\ot b\tilde]= (x\circ s)\ot\frac{1}{2}[a,b]+[x,s]\ot\frac{1}{2}(a\circ b)=-[s\ot b,x\ot a],\vspace{1mm}\\
\;[s\ot b,t\ot b'\tilde]=[s,t]\ot\frac{1}{2}(b\circ b')+ (s\circ t)\ot\frac{1}{2}[b,b']+tr(st)\la b,b'\ra_c,\vspace{1mm}\\
\;[x\ot a,u\ot c\tilde]=xu\ot a\cdot c=-[u\ot c,x\ot a],\vspace{1mm}\\
\;[s\ot b,u\ot c\tilde]=su\ot b\cdot c=-[u\ot c,s\ot b],\vspace{1mm}\\
\;[u\ot c,v\ot c'\tilde]=(u\circ v)\ot (c\diamond c')+ [u, v]\ot (c\heart c')+(u,v)\la c,c'\ra_c,\vspace{1mm}\\
\;[\la \b_1,\b_2\ra,x\ot a\tilde]=
\frac{-1}{4\ell}(x\circ
Id_{_{\v^\ell}}\ot[a,\b_{_{\b_1,\b_2}}^*]+[x,Id_{_{\v^\ell}}]\ot a\circ \b_{_{\b_1,\b_2}}^*),\vspace{1mm}\\
\;[\la \b_1,\b_2\ra,\hspace{-1mm}s\ot
b\tilde]\hspace{-1mm}=\hspace{-1mm}\frac{-1}{4\ell}([s,Id_{_{\v^\ell}}\hspace{-.5mm}]\hspace{-1mm}\ot\hspace{-1mm} (b\circ \b_{_{\b_1,\b_2}}^*\hspace{-1mm})\hspace{-1mm}+\hspace{-1mm}(s\circ\hspace{-1mm}
Id_{_{\v^\ell}})\hspace{-1mm}\ot \hspace{-1mm}[b, \b_{_{\b_1,\b_2}}^*\hspace{-1mm}]\hspace{-1mm}+\hspace{-1mm}2tr(sId_{\v^\ell})\la b,\b_{_{\b_1,\b_2}}^*\hspace{-.5mm}\ra_c),\vspace{1mm}\\
\;[\la \b_1,\b_2\ra_c,v\ot
c\tilde]=\frac{1}{2\ell}Id_{_{\v^\ell}}v\ot (\b_{_{\b_1,\b_2}}^*\cdot c)-\frac{1}{2}v\ot
(f(c,\b^*_2)\cdot \b^*_1+f(c,\b^*_1)\cdot \b^*_2)\\
\;[\la\b_1,\b_2\ra_c,\la\b'_1,\b'_2\ra_c\tilde]=\la d^\ell_{\b_1,\b_2}(\b'_1),\b'_2\ra_c+\la\b'_1,d^\ell_{\b_1,\b_2}(\b'_2)\ra_c
\end{array}
 \end{equation}
for $x,y\in\gg,$  $s,t\in\ss,$ $u,v\in\v,$
$a,a'\in\aa,$ $b,b'\in\bb,$ $c,c'\in\cc,$ $\b_1,\b_2,\b_1',\b'_2\in\fb.$
Moreover, under the above identification,   $\pi:\tilde\LL\longrightarrow \LL$ is given   by $\pi(x)=x$  for $x\in (\gg\ot \aa)\op(\ss\ot \bb)\op(\v\ot \cc)$ and   $\pi(\la b,b'\ra)=\la b,b'\ra_c$ for $b,b'\in \fb.$
\end{Theorem}
\pf As we have already seen, without loss of generality, we may assume $\tau(\gg,\LL)=\{0\}.$ We now note that $\tilde\LL$ is an $R-$graded  Lie algebra with grading pair $(\gg,\hh)$ and weight space decomposition  $\tilde\LL=\op_{\a\in R}\tilde\LL_\a$ where
\begin{equation}
\tilde\LL_\a=\LL_\a;\;\;\a\in R\setminus\{0\},\; \tilde\LL_0=\LL_0\op ker(\pi)=\LL_0\op E.
\end{equation}
 Suppose that $\{a_i\mid i\in\ii\},$ $\{b_j\mid j\in\jj\}$ and $\{c_t\mid t\in\T\}$ are bases for $\aa,\bb$ and $\cc$ respectively. We assume $0\in\ii$ and $a_0=1.$    For $\lam\in\Lam$ and $i\in\ii,j\in\jj$ and $t\in\T,$ we set
\begin{eqnarray*}
\gg^\lam_i:=\gg^\lam\ot a_i,\;\;\gg_i:=\gg\ot a_i\\
\ss^\lam_j:=\ss^\lam\ot b_j,\;\;\ss_j:=\ss\ot b_j\\
\v^\lam_t:=\v^\lam\ot c_t,\;\;\v_t:=\v\ot c_t.
\end{eqnarray*}
Therefore for $\dd:=\la\fb,\fb\ra,$ we have 
$\LL=\sum_{i\in\i}\gg_i\op\sum_{j\in \jj}\ss_j\op\sum_{t\in\T}\v_t\op\dd,$  $\gg_i=\cup_{\lam\in\Lam}\gg_i^\lam,$ $\ss_j=\cup_{\lam\in\Lam}\ss_j^\lam$ and $\v_t=\cup_{\lam\in\Lam}\v_t^\lam,$ $i\in\ii,j\in\jj,t\in\T.$ For $\lam\in\lam,$  set $$\begin{array}{l}\LL^\lam:=\sum_{i\in\i}\gg^\lam_i\op\sum_{j\in \jj}\ss^\lam_j\op\sum_{t\in\T}\v^\lam_t\op\dd,\\\\
\hat\LL^\lam:= \sum_{i\in\i}\gg^\lam_i\op\sum_{j\in \jj}\ss^\lam_j\op\sum_{t\in\T}\v^\lam_t\op\dd\op ker(\pi),\\\\
\tilde\LL^\lam:=[\hat\LL^\lam,\hat\LL^\lam\tilde]=[\LL^\lam,\LL^\lam\tilde].\end{array}$$ The restriction of $\pi$ to $\hat\LL^\lam$ is a central extension of $\LL^\lam $ and setting   $\pi_\lam:=\pi|_{_{\tilde\LL^\lam}}:\tilde\LL^\lam\longrightarrow \LL^\lam,$ we get that  $(\tilde\LL^\lam,\pi_\lam)$ is a perfect central extension of $\LL^\lam.$ Also by Lemma \ref{uce8}, we have \begin{equation}\label{new6}\tilde \LL^\lam=\LL^\lam\op \z_\lam\end{equation} where $\z_\lam:=ker(\pi_\lam).$
Now as  $\tilde\LL$ is perfect,
\begin{eqnarray*}\tilde\LL=[\tilde\LL,\tilde\LL\tilde]=[\cup_{\lam\in\Lam}\hat\LL^\lam,\cup_{\lam\in\Lam}\hat\LL^\lam]
=\cup_{\lam\in\Lam}[\hat\LL^\lam,\hat\LL^\lam]
=\cup_{\lam\in\Lam}\tilde\LL^\lam
\end{eqnarray*}
so  $\tilde\LL$ is the direct union of $\{\tilde\LL^\lam\mid \lam\in\Lam\}.$ We next note that $\LL^\lam$ is an $R_\lam-$graded Lie algebra with grading pair $(\gg^\lam,\hh_\lam:=\gg^\lam\cap\hh)$  and $\tilde\LL^\lam$ is a perfect central extension of $\LL^\lam$ with corresponding $2-$cocycle $\tau_\lam:=\tau\mid_{_{\LL^\lam\times\LL^\lam}}$ satisfying $\tau_\lam(\gg^\lam,\LL^\lam)=\{0\}.$ Therefore by Lemma \ref{first}, $\tilde\LL^\lam=\op_{\a\in R_\lam}\tilde\LL^\lam_\a$ with \begin{equation*}
\tilde\LL^\lam_\a=\left\{\begin{array}{ll}\LL^\lam_\a=\LL_\a& \hbox{if $\a\in R_\lam\setminus\{0\}$}\\
\LL^\lam_0\op\z_\lam& \hbox{if $\a=0$}
\end{array}
\right.=\left\{\begin{array}{ll}\LL_\a& \hbox{if $\a\in R_\lam\setminus\{0\}$}\\
\displaystyle{\sum_{\b\in R_\lam\setminus\{0\}}}[\LL_\b,\LL_{-\b}\tilde]&\hbox{if $\a=0.$}
\end{array}
\right.
\end{equation*}
We next recall from (\ref{new6}) that for $\lam\in\Lam,$ $\tilde\LL^\lam=\sum_{i\in \i}\gg_i^\lam\op\sum_{j\in \jj}\ss_j^\lam\op\sum_{t\in \T}\v_t^\lam\op \dd_\lam\op\z_\lam.$ Also  $\dd_\lam\op\z_\lam$ is a trivial $\gg^\lam-$submodule of $\tilde\LL^\lam.$  We next note that $\tau(\gg,\LL)=\{0\}$ which implies that $\LL^\lam$ is a $\gg^\lam-$submodule of $\tilde\LL^\lam.$ Now as  $\gg_i^\lam$ is the $\gg^\lam-$submodule of $\LL^\lam$ generated by $\gg_i^0,$ we get that $\gg_i^\lam$ is also the  $\gg^\lam-$submodule  of $\tilde\LL^\lam$ generated by $\gg_i^0.$ Similarly for $j\in\jj$ and $t\in\T,$ $\ss_j^\lam$ coincides with the $\gg^\lam-$submodule of  $\tilde\LL^\lam$ generated by $\ss_j^0,$ and $\v_t^\lam$ coincides with the $\gg^\lam-$submodule of  $\tilde\LL^\lam$ generated by $\v_t^0.$
This means that $$(\i,\jj,\T,\{\gg_i^0\},\{\gg_i^\lam\},\{\ss_j^0\},\{\ss_j^\lam\},\{\v_t^0\},\{\v_t^\lam\},\dd_0\op\z_0,\dd_\lam\op\z_\lam)$$ is an $(R_0,R_\lam)-$datum for $0\prec\lam$ in the sense of \cite{Y}. Therefore using [Y], we get that
\begin{eqnarray*}\tilde\LL^\lam&=&\sum_{i\in \i}\gg_i^\lam\op\sum_{j\in \jj}\ss_j^\lam\op\sum_{t\in \T}\v_t^\lam\op \dd_0\op\z_0\\
&=&(\gg^\lam\dot\ot\aa)\op(\ss^\lam\dot \ot \bb)\op(\v^\lam\dot \ot \cc)\op\dd_0\op\z_0.
\end{eqnarray*} So
\begin{eqnarray*}
\tilde\LL=\cup\tilde\LL^\lam&=&\sum_{i\in \i}\gg_i\op\sum_{j\in \jj}\ss_j\op\sum_{t\in \T}\v_t\op \dd_0\op\z_0\\
&=&(\gg\ot\aa)\op(\ss\ot \bb)\op(\v\ot \cc)\op\dd_0\op\z_0.
\end{eqnarray*}
Using [Y] and Lemma \ref{first}, $\LL^0,$ $\tilde\LL^0,$ $\tilde\LL^\lam$ and $\LL^\lam$ have the same coordinate quadruple  $\mathfrak{q}.$ Also there is a subspace $\kk_0$ of ${\rm HF}(\fb)$ satisfying the uniform property on $\fb$ such that $\dd_0\op\z_0=\{\fb,\fb\}/\kk_0.$ Now setting $\la\fb,\fb\ra_c:=\{\fb,\fb\}/\kk_0,$ we get using  (\ref{probc-fin}) using [Y, Pro 3.10]. Now by [Y, Theorem 4.1], $\tilde\LL$ is an $R-$graded Lie algebra.
Now for fix $x,y\in\gg$ with $tr(xy)\neq 0,$ and $a,a'\in\aa,$ we have
\begin{eqnarray*}
&&[x,y]\ot(1/2)(a\circ a')+(x\circ y)\ot (1/2)[a,a']+tr(xy)\la a,a'\ra\\
&=&[x\ot a,y\ot a']\\&=&[\pi(x\ot a),\pi(y\ot a')]\\
&=&\pi([x\ot a,y\ot a'\tilde])\\
&=&\pi([x,y]\ot(1/2)(a\circ a')+(x\circ y)\ot (1/2)[a,a']+tr(xy)\la a,a'\ra_c)\\&=&
[x,y]\ot(1/2)(a\circ a')+(x\circ y)\ot (1/2)[a,a']+tr(xy)\pi(\la a,a'\ra_c).
\end{eqnarray*}
This implies that $\pi(\la a,a'\ra_c)=\la a,a'\ra.$ Similarly we can prove that $\pi(\la b,b'\ra_c)=\la b,b'\ra$ and $\pi(\la c,c'\ra_c)=\la c,c'\ra$ for $b,b'\in\bb$ and $c,c'\in\cc.$ This completes the proof.\qed

\begin{Theorem} Suppose that $\mathfrak{q}:=(\fa,*,\cc,f)$ is a coordinate quadruple of type $BC,$ $\fb:=\fb(\mathfrak{q}),$ $\kk$ a subspace of ${\rm HF}(\fb)$ satisfying the uniform property on $\fb$ and   $\LL:=\LL(\mathfrak{q},\kk)=(\gg\ot\aa)\op(\ss\ot\bb)\op(\v\ot\cc)\op\la\fb,\fb\ra,$ where $\la\fb,\fb\ra:=\{\fb,\fb\}/\kk,$  the corresponding  $R-$graded Lie algebra. Consider Remark \ref{rem1} and set $\mathfrak{A}:=\LL(\mathfrak{q},\
\{0\})=(\gg\ot\aa)\op(\ss\ot\bb)\op(\v\ot\cc)\op\{\fb,\fb\},$ then $\mathfrak{A}$ is the universal central extension of $\LL.$
\end{Theorem}
\pf  Define \begin{equation*}
\begin{array}{l}\pi:\mathfrak{A}\longrightarrow \LL;\\
x\mapsto x;\;\; x\in (\gg\ot\aa)\op(\ss\ot\bb)\op(\v\ot\cc);\\
\{ b,b'\}\ra_u\mapsto \{b,b'\}+\kk=\la b,b'\ra.\end{array}
\end{equation*}
If $x\in ker(\pi),$ then $x=\sum\{ \b_i,\b'_i\}$ such that $\sum_i\la \b_i,\b'_i\ra=0.$ But since $\kk$ satisfies the uniform property on $\fb,$ we get that $\sum_i\b^*_{\b_i,\b'_i}=0.$ Now (\ref{probc-fin}) together with (\ref{probc-gen}) implies that $x\in Z(\mathfrak{A}).$ This means that $\pi$ is a central extension of $\LL.$ Now suppose that $\dot\LL$ is a  Lie algebra and $\dot\varphi:\dot\LL\longrightarrow \LL$ is a central extension of $\LL.$ Set $\tilde\LL$ to be the derived algebra of $\dot\LL$ and $\varphi:=\dot\varphi\mid_{_{\tilde\LL}}.$ Then $(\tilde\LL,\varphi)$ is a perfect central extension of $\LL.$ By Theorem \ref{main'}, we may assume there is a subspace $\kk_0$ of  ${\rm HF}(\fb)$ satisfying the uniform property on $\fb$ such that $\tilde\LL=(\gg\ot\aa)\op(\ss\ot\bb)\op(\v\ot\cc)\op\la\fb,\fb\ra_c$ where $\la\fb,\fb\ra_c:=\{\fb,\fb\}/\kk_0$ and $\varphi:\tilde\LL\longrightarrow \LL$ is given by $$\begin{array}{l}\varphi(x)=x;\;\; x\in (\gg\ot\aa)\op(\ss\ot\bb)\op(\v\ot\cc)\\
\varphi(\la \b,\b'\ra_c)=\la \b,\b'\ra;\;\b,\b'\in\fb.\end{array}$$ Now if we define
$$\begin{array}{l}\psi:(\gg\ot\aa)\op(\ss\ot\bb)\op(\v\ot\cc)\op\{\fb,\fb\}\longrightarrow \tilde\LL\\
\psi(x)=x;\;\;x\in (\gg\ot\aa)\op(\ss\ot\bb)\op(\v\ot\cc)\andd \psi(\{ \b,\b'\})=\la \b,\b'\ra_c,
\end{array}$$
$\psi$ is a Lie algebra homomorphism satisfying $\varphi\circ\psi=\pi.$ In other words, $\pi$ is the universal central extension.\qed

\vspace{2cm}

 Department of Mathematics, University of Isfahan, Isfahan, Iran,
P.O.Box 81745-163.\\            
ma.yousofzadeh@sci.ui.ac.ir               
\\
 School of Mathematics, Institute for Research in
Fundamental
Sciences (IPM), P.O. Box: 19395-5746, Tehran, Iran.\\

\end{document}